\documentclass{emsprocart}

\usepackage{diagrams,graphicx,color}
\usepackage[color,matrix,arrow]{xy}
\usepackage{makeidx}

\usepackage[metapost,truebbox]{mfpic}
\opengraphsfile{mypicts}

\xyoption{curve}
\xyoption{arrow}
\xyoption{matrix}

\contact[lieven.lebruyn@ua.ac.be]{Lieven Le Bruyn \\ Department of Mathematics, University of Antwerp \\ Middelheimlaan 1, B-2012 Antwerp, Belgium}

\newcommand{\wis}[1]{{\text{\em \usefont{OT1}{cmtt}{m}{n} #1}}}
\newcommand{\V}{\mathbb{V}}
\newcommand{\N}{\mathbb{N}}
\newcommand{\C}{\mathbb{C}}
\newcommand{\Z}{\mathbb{Z}}
\newcommand{\F}{\mathbb{F}}
\newcommand{\Q}{\mathbb{Q}}
\newcommand{\PP}{\mathbb{P}}

\newcommand{\Oscr}{\mathcal{O}}

\theoremstyle{definition}

\makeindex

\title[Absolute geometry and the Habiro topology]{Absolute geometry and the Habiro topology}

\author[Lieven Le Bruyn]{Lieven Le Bruyn}
 
\begin{document}

\begin{quote}
"One can postulate, of course, that $\wis{Spec}(\F_1)$ is the absolute point, but the real problem is to develop non-trivial consequences of this point of view."

Mikhail Kapranov and Alexander Smirnov in \cite{SmirnovKapranov96}
\end{quote}

\begin{quote}
"Analogies, it is true, decide nothing, but they can make one feel more at home."

Sigmund Freud, The Essentials of Psycho-Analysis
\end{quote}

\vskip 4mm
Ever since the work of Richard Dedekind and Leopold Kronecker the striking similarities between number fields and function fields of smooth projective curves over finite fields have served as a powerful analogy to transport results and conjectures from number theory to these function fields. Using the powerful machinery of geometry one could prove results on the function field side for which the corresponding number theoretic result still remains conjectural. Most noteworthy are Andre Weil's proof of the Riemann hypothesis for function fields in the 40ties and, more recently, the proof of the ABC-conjecture for function fields. 

Since the mid 80ties, attempts have been made to mimic Weil's approach to the Riemann hypothesis by trying to imagine the integral prime spectrum $\wis{Spec}(\Z)$ to be a 'curve' over the {\em absolute point} \index{absolute point}, that is $\wis{Spec}(\F_1)$ where $\F_1$ is the elusive field with one element, and subsequently to study intersection theory on the 'surface' $\wis{Spec}(\Z) \times_{F_1} \wis{Spec}(\Z)$, which is part of Weil's approach to the Riemann hypothesis. Since he could not invent what this 'surface' might be, Alexander Smirnov decided to study intersection theory on the easier 'surface' $\PP^1_{\F_1} \times_{\F_1} \wis{Spec}(\Z)$, and in particular to investigate the graphs of 'maps'
\[
q~:~\wis{Spec}(\Z) \rTo \PP^1_{\F_1} \qquad \text{where $q \in \Q$} \]
which should exist by analogy with the function field case as $\Q$ should be thought of as the function field of the 'curve' $\wis{Spec}(\Z)$. Smirnov dreamed up the following definitions for both geometric objects and of the maps between them (see \cite{Smirnov92}):
\begin{itemize}
\item{The {\em absolute projective line} \index{absolute projective line}, $\PP^1_{\F_1}$, should have as its schematic points the set
\[
\{ [0],[\infty] \} \cup \{ [n]~|~n \in \N_0 \} \]
where the degree of the point $[n]$ should be $\phi(n)$ for $n \in \N_0$ and equal to $1$ for $[0]$ and $[\infty]$. The schematic point $[n]$ should be thought of as corresponding to the set of geometric points being the primitive $n$-th roots of unity.}
\item{The {\em completed prime spectrum} \index{completed prime spectrum}, $\overline{\wis{Spec}(\Z)}$ should have as its schematic points the set
\[
\{ (p)~|~p~\text{a prime number} \} \cup \{ \infty \} \]
where the degree of $\infty$ is equal to one and the degree of $(p)$ should be $log(p)$.}
\item{If $q = \tfrac{a}{b} \in \Q$ with $(a,b)=1$, the map
\[
q~:~\overline{\wis{Spec}(\Z)} \rTo \PP^1_{\F_1} \]
should send the point $(p)$ to $[0]$ if $p$ is a prime factor of $a$, to $[\infty]$ if $p$ is a prime factor of $b$ and, in the remaining cases, to $[n]$ if $n$ is the order of the image of $\tfrac{a}{b}$ in the finite group $\F_p^*$. Finally, $\infty$ should be send to $[0]$ if $a<b$ and to $[\infty]$ if $a > b$.}
\end{itemize}
Here's part of the graph of the map $q=2$ (for primes $p < 1000$) in the 'surface' $\PP^1_{\F_1} \times \overline{\wis{Spec}(\Z)}$
\[
\includegraphics[width=12cm]{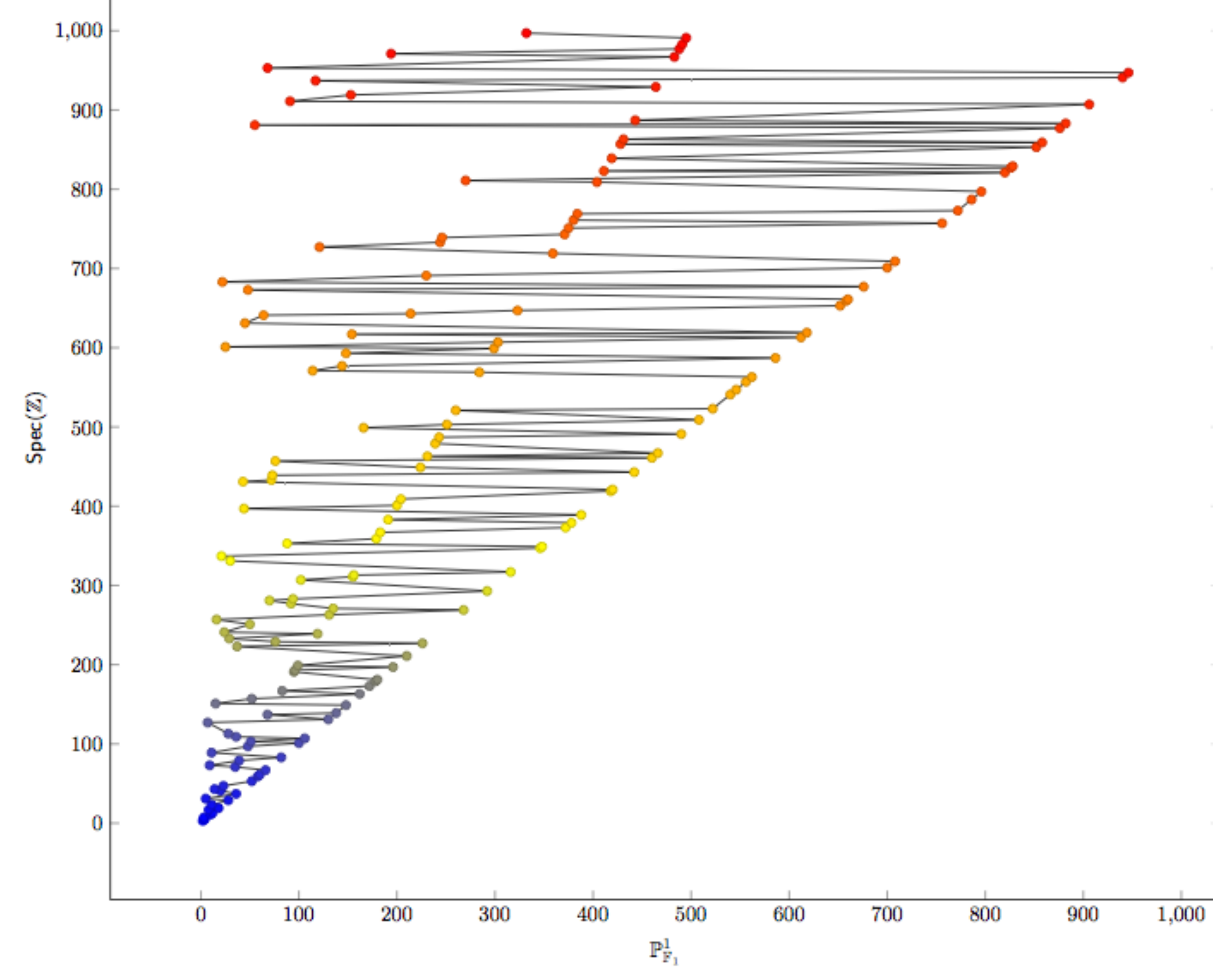}
\]
It is easy to verify that these maps are finite but showing that they are actually covers for most values of $q$ relies upon a result of Zsigmondy's \cite{Zsig}. In \cite{Smirnov92} Alexander Smirnov was able to deduce the ABC-conjecture for $\Z$ provided one would be able to develop an absolute geometry, admitting suitable versions of $\wis{Spec}(\Z)$ and $\PP^1_{\F_1}$ and such that one can prove an analogue of the Riemann-Hurwitz formula for maps such as $q$.
Since then, numerous proposals for a geometry over $\F_1$ have been made, all of them allowing objects such as $\PP^1_{\F_1}$ and similar combinatorial defined varieties such as affine and projective spaces, Grassmannians etc., but almost none of them containing objects having the desired properties of $\overline{\wis{Spec}(\Z)}$. For an overview of these attempts and the connections between them we refer to \cite{PenaLorscheid09}. Perhaps the most promising approach was put forward by Jim Borger and is based on the notion of $\lambda$-rings, \cite{Borger09}. For our purposes, a $\lambda$-ring is a $\Z$-algebra $R$ without additive torsion admitting a commuting family of endomorphisms $\{ \Psi^n~:~n \in \N_0 \}$ such that for prime numbers $p$ the map $\Psi^p$ is a lift of the Frobenius morphism on $R \otimes_{\Z} \F_p$. Borger interprets this family of endomorphisms as descent data from $\Z$ to $\F_1$, and conversely views the forgetful functor, stripping off the $\lambda$-structure, as base-extension $- \otimes_{\F_1} \Z$. In this approach, $\PP^1_{\F_1}$ would then be the usual integral scheme $\PP^1_{\Z}$ equipped with the toric $\lambda$-ring structure induced by the endomorphisms $\Psi^n(x)=x^n$ on $\Z[x]$, giving us the fanciful identity
\[
\PP^1_{\F_1} \times_{\wis{Spec}(\F_1)} \wis{Spec}(\Z) = \PP^1_{\Z} \]
where on the right-hand side we forget the $\lambda$-structure on the integral projective line $\PP^1_{\Z}$, giving us a concrete proposal for Smirnov's plane. Unlike other approaches, Borger's proposal allows us define how an integral scheme $X_{\Z}$ should be viewed over $\F_1$. Indeed, the forgetful-functor (that is, base-extension $- \otimes_{\F_1} \Z$) has a right-adjoint functor 
\[
w~:~\wis{rings} \rTo \wis{rings}_{\lambda} \]
assigning to a $\Z$-algebra $A$ a close relative of the ring of Big Witt vectors, $w(A) = 1 + tA[[t]]$, equipped with a new addition $\oplus$ being the ordinary multiplication of power series, and a new multiplication $\otimes$ induced functorially by the condition that
\[
(\frac{1}{1-a.t}) \otimes (\frac{1}{1-b.t}) = \frac{1}{1-ab.t} \]
for all $a,b \in A$. This functor can then be viewed as Weil-restriction from integral schemes to $\F_1$-schemes. Hence, in particular, this proposal allows us to define $\wis{Spec}(\Z)/\F_1$ as the $\F_1$-geometric object corresponding to the ring $w(\Z)$ which is isomorphic to the completed Burnside-ring $\hat{B}(C)$ of the infinite cyclic group $C$, by \cite{DressSieb89}. In these notes we will explore how Smirnov's maps $\wis{Spec}(\Z) \rOnto \PP^1_{\F_1}$ fit into Borger's proposal.

A second theme of these notes is to explore the origins of a new topology on the roots of unity $\pmb{\mu}_{\infty}$ introduced and studied by Kazuo Habiro in \cite{Habiro02} in order to unify invariants of $3$-dimensional homology spheres, introduced first by Ed Witten by means of path integrals and rigorously constructed by Reshitikhin and Turaev.
Habiro calls two roots of unity adjacent to each other whenever their quotient is of pure prime power order. For example below we depict the adjacency relation on $60$-th roots of unity where we used different colors for different prime powers ($2$-powers are colored yellow, $3$- and $5$-powers respectively blue and red).
\[
\includegraphics[height=7cm]{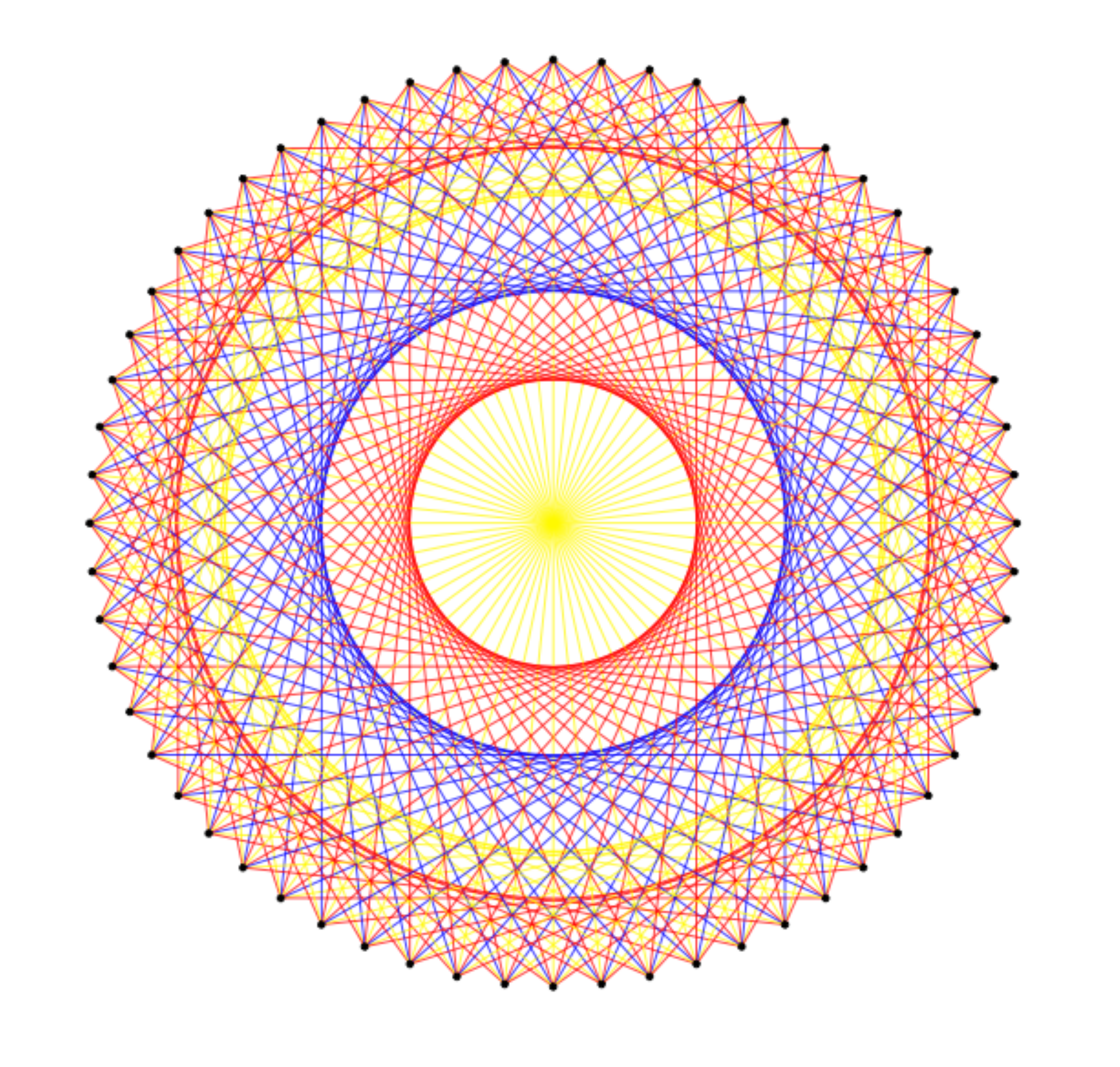} \]
The {\em Habiro topology} on $\pmb{\mu}_{\infty}$ is then defined by taking as its open sets those subsets $U \subset \pmb{\mu}_{\infty}$ having the property that for every $\alpha \in U$ all except finitely many $\beta \in \pmb{\mu}_{\infty}$ adjacent to $\alpha$ also belong to $U$. The Galois action is continuous in this topology which is in sharp contrast to the induced analytic topology. The Habiro topology is best understood by applying techniques from noncommutative algebraic geometry to objects like $\PP^1_{\F_1}$. Recall that the schematic point $[n]$ of $\PP^1_{\F_1}$ corresponds to the set of primitive $n$-th roots of unity and hence corresponds to the closed subscheme of $\PP^1_{\Z}$ defined by the $n$-th cyclotomic polynomial $\Phi_n(x)$. For $n \not= m$ the corresponding ideals do not have to be co-maximal (that is, the closed subschemes can intersect over some prime numbers $p$) and, in fact, whenever $\tfrac{m}{n}=p^k$ for some prime number $p$ there are non-split extensions of $\Z[x,x^{-1}]$-modules
\[
0 \rTo \Z[\zeta_n] \rTo E \rTo \Z[\zeta_m] \rTo 0 \]
In noncommutative algebraic geometry such situations are interpreted as saying that the corresponding points $[m]$ and $[n]$ lie infinitely close to each other as they share some common tangent information. This then is the origin of the Habiro topology on $\pmb{\mu}_{\infty}$. So in the above picture one should view two roots of unity to be infinitely close whenever they are connected by a colored line, giving us a horrible topological space. The tools of noncommutative geometry allow us to study such bad spaces by associating noncommutative algebras to them and in this case, the Bost-Connes algebra $\Lambda$ is naturally arises from it. More generally, one assigns to a $\lambda$-ring $A$ a noncommutative algebra, the skew-monoid algebra $A \ast \N^{\times}_0$ where the skew-action is determined by the family of endomorphisms $\Psi^n$. Therefore, one might argue that $\F_1$-geometry is essentially of a  noncommutative nature. In these notes we will explore this line of thoughts and show in particular that the Habiro topology on $\PP^1_{\Z}$ is a proper refinement of the Zariski- (that is, cofinite-) topology and is no longer compact. We can then also define an exotic new topology on $\wis{Spec}(\Z)$ by demanding that all the Smirnov-maps $q~:~\wis{Spec}(\Z) \rTo \PP^1_{\F_1}$ should be continuous with respect to the Habiro topology on $\PP^1_{\F_1}$.

\tableofcontents

\par \vfill \eject

\section{Mumford's drawings of $\mathbb{A}^1_{\Z}$ and $\PP^1_{\Z}$}

Let us start with the iconic drawing of the 'arithmetic surface', that is of the prime spectrum $\mathbb{A}^1_{\Z} = \wis{Spec}(\Z[x])$, by David Mumford in the original version of his Red Book \cite[p. 141]{Mumf67}. Subsequent more polished versions of the drawing can be found in the reprinted Red Book \cite[p. 75]{Mumf99} and in a.o. \cite[p.24]{Reid95} and \cite[p. 85]{EisenbudHarris00}. 

\[
\includegraphics[width=12cm]{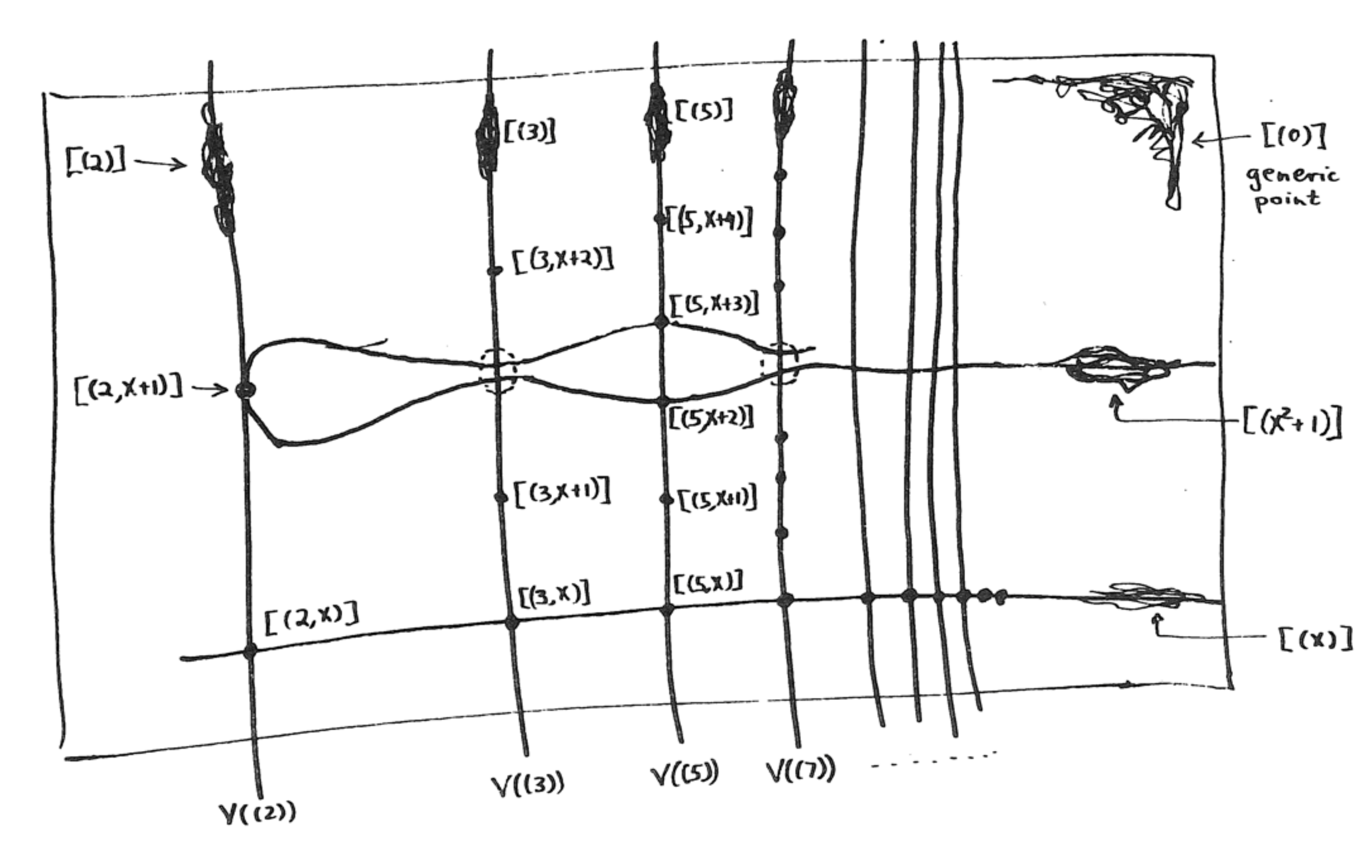} \]
It was believed to be the first depiction of one of Grothendieck's prime spectra having a real mixing of arithmetic and geometric properties, and as such was influential for generations of arithmetic geometers. Clearly, $\mathbb{A}^1_{\Z} = \wis{Spec}(\Z[x])$, that is the set of all prime ideals of $\Z[x]$, contains the following elements
\begin{itemize}
\item{$(0)$ depicted as the generic point $[(0)]$}
\item{principal prime ideals $(f)$, where $f$ is either a prime number $p$ (giving the vertical lines $\mathbb{V}((p)) = \wis{Spec}(\mathbb{F}_p[x])$) or a $\mathbb{Q}$-irreducible polynomial written so that its coefficients have greatest common divisor $1$ (the horizontal 'curves' in the picture such as $[(x^2+1)]$}
\item{maximal ideals $(p,f)$ where $p$ is a prime number and $f$ is a monic polynomial which remains irreducible modulo $p$, the 'points' in the picture.}
\end{itemize}
Mumford's drawing focusses on the vertical direction as the vertical lines $\mathbb{V}((p))$ are the fibers of the projection $\wis{Spec}(\mathbb{Z}[x]) \rOnto \wis{Spec}(\mathbb{Z})$ associated to the structural map $\Z \rInto \Z[x]$. This is consistent with Mumford's drawing of $\wis{Spec}(\Z)$ in \cite[p. 137]{Mumf67} where he writes "$\Z$ is a principal ideal domain like $k[x]$, and $\wis{Spec}(\Z)$ is {\em usually} visualized as a line :
\[
\includegraphics[width=12cm]{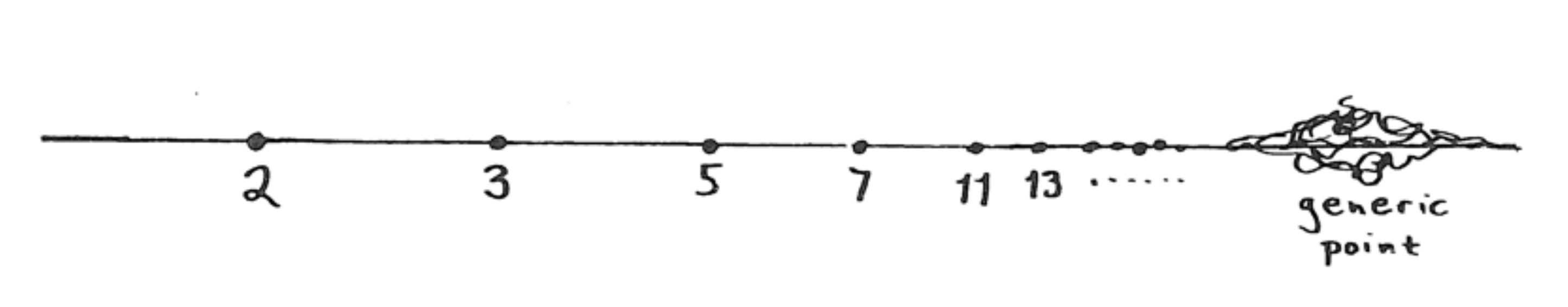} \]
There is one closed point for each prime number, plus a generic point $[(0)]$." I've emphasized the word 'usually' as Mumford knew at the time he was writing the Red Book perfectly well that there were other, and potentially better, descriptions of $\wis{Spec}(\Z)$ than this archaic prime number line. 

In july 1964 David Mumford attended the Woods-Hole conference, which became famous for producing the Atiyah-Bott fixed point theorem.
On july 10th there were three talks on the hot topic at the time emerging from Grothendieck's Parisian seminar : Etale cohomology. Mike Artin spoke on 'Etale cohomology of schemes' \cite{Artin64}, Jean-Louis Verdier on 'A duality theorem in the etale cohomology of schemes' \cite{Verdier64} and John Tate on 'Etale cohomology over number fields' \cite{Tate64}. Later in the conference, Mike Artin and Jean-Louis Verdier ran a 'Seminar on \'etale cohomology of number fields' \cite{ArtinVerdier64} in which they proved their famous duality result
\[
H^r_{et}(\wis{Spec}(\Z),\mathcal{F}) \times Ext^{3-r}_{\wis{Spec}(\Z)}(\mathcal{F},\mathbb{G}_m) \rTo H^3_{et}(\wis{Spec}(\Z),\mathcal{F}) \simeq \Q/\Z \]
for Abelian constructible sheaves $\mathcal{F}$, suggesting a $3$-dimensional picture of $\wis{Spec}(\Z)$. 

Combining this with the fact that the \'etale fundamental group of $\wis{Spec}(\Z)$ is trivial (and that of $\wis{Spec}(\mathbb{F}_p)$ is the profinite completion of $\pi_1(S^1)=\Z$), Mumford dreamed-up the analogy between prime number and knots in $S^3$, see for example the opening paragraph of the unpublished preprint by Barry Mazur 'Remarks on the Alexander polynomial' \cite{Mazur65}:
"Guided by the results of Artin and Tate applied to the calculation of the Grothendieck Cohomology Groups of the schemes
\[
\wis{Spec}(\Z/ p \Z) \rInto \wis{Spec}(\Z) \]
Mumford has suggested a most elegant model as a geometric interpretation of the above situation : $\wis{Spec}(\Z/p\Z)$ is like a one-dimensional knot in $\wis{Spec}(\Z)$ which is like a simply connected three-manifold." This analogy between prime numbers and knots has led in the past decades to the field of 'Arithmetic Topology', a good introduction to which can be found in the lecture notes by Masanori Miyashita \cite{Miyashita09}. 

However, the arithmetic plane wasn't the first attempt by Mumford to draw an arithmetic scheme. In his Harvard 'Lectures on curves on an algebraic surface' \cite{Mumf66} there's on page 28 this drawing of $\PP^1_{\Z}=\wis{Proj}(\Z[X,Y])$
\[
\includegraphics[width=12cm]{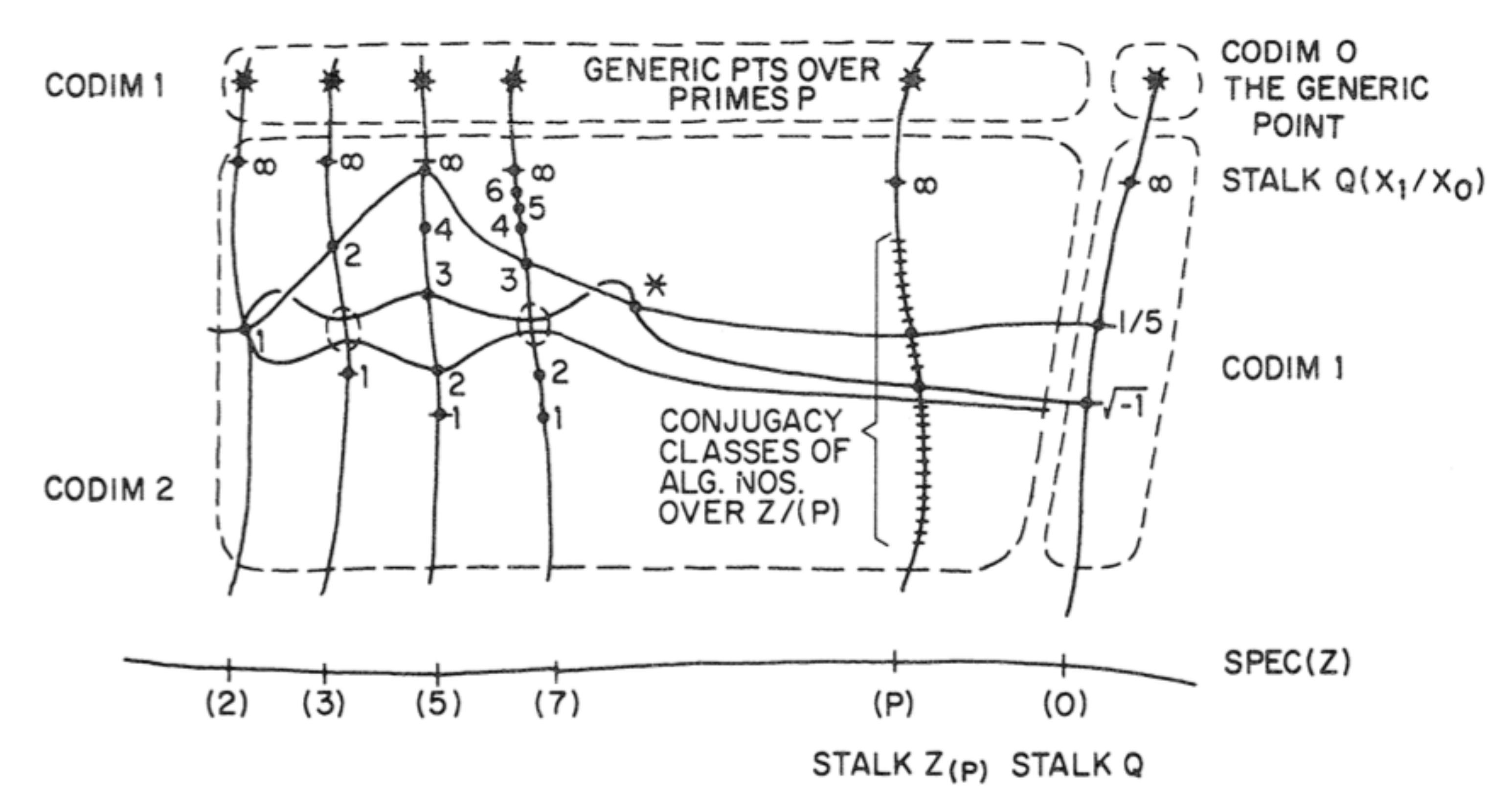}
\]
This drawing has at the same time a more classical touch to it, separating the different elements of $\wis{Proj}(\Z[X,Y])$ (that is the graded prime ideals of $\Z[X,Y]$ not containing in $(X,Y)$) according to codimension, as well as being more modern in that there is a $3$-dimensional feel to it (the closed subschemes $\mathbb{V}(X^2+Y^2)$ and $\mathbb{V}(5X-Y)$ have over- and undercrossings). The points of $\wis{Proj}(\Z[X,Y])$ are
\begin{itemize}
\item{The graded ideal $0$ corresponding to the unique codimension zero point, the generic point,}
\item{The graded prime ideals of height one which are either the vertical fibers $\mathbb{V}((p)) = \PP^1_{\F_p} = \wis{Proj}(\F_p[X,Y])$ or the horizontal subschemes corresponding to the homogenization (with respect to $Y$) of a $\Q$-irreducible polynomial in $\Z[X]$ written such that the greatest common divisor of its coefficients equals $1$, the codimension one points, and,}
\item{The codimension two points, which correspond to the graded ideals $(p,F)$ where $F$ is a homogeneous element of $\Z[X,Y]$ such that its reduction modulo $p$ remains irreducible.}
\end{itemize}
For example, the point marked $\ast$ in the drawing is the point $(13,X-8Y)$.
This picture resembles that of $\mathbb{A}^1_{\Z}$ and is in fact the gluing of two such drawings over their intersection, the first is obtained by removing the $\infty$-section (that is $\mathbb{V}(Y)$) and is $\wis{Spec}(\Z[x])$ with $x = \tfrac{X}{Y}$ whereas the second is obtained by removing the $0$-section $\mathbb{V}(X)$ and is $\wis{Spec}(\Z[x^{-1}])$. They are glued together over their intersection $\wis{Spec}(\Z[x,ux^{-1}])$. 

Influential as these drawings have been, there are a couple of obvious problems with them which will lead us unavoidably to the concept of the {\em absolute point} $\wis{Spec}(\F_1)$, that is, the geometric object associated to the elusive {\em field with one element} $\F_1$.

\begin{enumerate}
\item{{\bf what is the vertical axis?} These drawings of $\mathbb{A}^1_{\Z}$ and $\PP^1_{\Z}$ as arithmetic 'planes' suggest that, apart from the 'horizontal axis' $\wis{Spec}(\Z)$, coming from the structural morphisms $\Z \rInto \Z[x^{\pm}]$, there should also be a 'vertical axis' and corresponding projection, so what is it?}
\item{{\bf what is the correct topology?} The drawing of the horizontal curves suggests a natural identification between vertical fibers $\PP^1_{\F_q} \leftrightarrow \PP^1_{\F_q}$ for primes $p \not= q$, so is there one? And what is the correct topology on these fibers, and on $\wis{Spec}(\Z)$?}
\end{enumerate}

\section{The vertical axis $\PP^1_{\F_1}$}

We have seen that the 'points' correspond to maximal ideals of $\Z[x]$ which are all of the form $\mathfrak{m}=(p,F)$ where $p$ is a prime number and $F$ is a monic irreducible polynomial such that its reduction $\overline{F} \in \mathbb{F}_p[x]$ remains irreducible. Clearly $\mathfrak{m}$ lies on a unique vertical ruling $\mathbb{V}((p))$ and we wonder whether there exists an appropriate set of horizontal rulings containing all points $\mathfrak{m}$. We know that the quotient
\[
\frac{\Z[x]}{\mathfrak{m}} \simeq \frac{\mathbb{F}_p[x]}{(\overline{F})} \simeq \mathbb{F}_{p^d} \]
is the finite field $\mathbb{F}_{p^d}$ where $d$ is the degree of $F$ and that its multiplicative group of units is the cyclic group $C_{p^d-1}$ and hence $\mathbb{F}_{p^n}$ consists of roots of unity together with zero.

This observation led Yuri I. Manin in \cite{Manin08} to consider the ring $\Z[x]_S$ which is the localization of $\Z[x]$ at the multiplicative system $S$ generated by the polynomials $\Phi_0(x)=x$ together with the {\em cyclotomic polynomials} \index{cyclotomic polynomial}
\[
\Phi_n(x) = \prod_{i=1}^{\mu(n)} (x-\epsilon_i) \]
where $\epsilon_i$ runs over all primitive roots of unity of order $n$, of which there are exactly $\phi(n)$ where $\phi$ is the Euler-function $\phi(n) = \# \{ 1 \leq j < n~:~(j,n)=1 \}$. It follows that the above point $\mathfrak{m}$ lies on the 'curve' $[(\Phi_{p^d-1}(x))]$ in the arithmetic plane $\wis{Spec}(\Z[x])$. Hence, localizing at $S$ removes all these curves $[(\Phi_n(x))]$ together with all the points lying on them. That is, $\wis{Spec}(\Z[x]_S)$ has no height two prime ideals and so consists of $(0)$ and the remaining height one prime ideals, all of which are principal. We conclude that the localized ring $\Z[x]_S$ is a principal ideal domain.

This is completely analogously to the more classical setting in which we localize $\Z[x]$ at the multiplicative system $S'$ generated by all prime numbers $p$, thus getting the principal ideal domain $\Q[x]$. 
Here the localization removes all the vertical rulings $\V((p))$ from the arithmetic plane together with all the points $\mathfrak{m}$ lying on them. 
Hence, this suggests to take the 'curves' $[(\Phi_n(x))] = \mathbb{V}((\Phi_n(x)))$ as a set of horizontal rulings, and then indeed the point $\mathfrak{m}$ lies on the intersection of the vertical ruling $\mathbb{V}((p))$ and the horizontal ruling $\mathbb{V}((\Phi_{p^d-1}(x)))$. 

Manin writes : "This suggests that the union of cyclotomic arithmetic curves $\Phi_n(x)=0$ can be imagined as the union of closed fibers of the projection
\[
\wis{Spec}(\Z[x]) \rOnto \wis{Spec}(\mathbb{F}_1[x]) \]
and the arithmetic plane itself as the product of two coordinate axes, an arithmetic one, $\wis{Spec}(\Z)$, and a geometric one, $\wis{Spec}(\mathbb{F}_1[x])$, over the 'absolute point' $\wis{Spec}(\mathbb{F}_1)$." 

Clearly we can repeat the same argument for $\wis{Spec}(\Z[x^{-1}])$ and we obtain as Manin's proposal for a set of horizontal rulings on $\PP^1_{\Z}$ the set of codimension one closed subschemes determined by the irreducible homogeneous polynomials
\[
\{ \Phi_0=X, \Phi_{\infty}= Y \} \cup \{ \Phi_n = \prod_{i=0}^{\mu(n)} (X-\epsilon_iY)~:~n \in \N_0 \} \]
That is, we can extend in Figure 1 Mumford's drawing of $\PP^1_{\Z}$ with an horizontal projection such that every codimension two point of $\PP^1_{\Z}$ lies on the intersection of a vertical and an horizontal ruling
\begin{figure}[t]
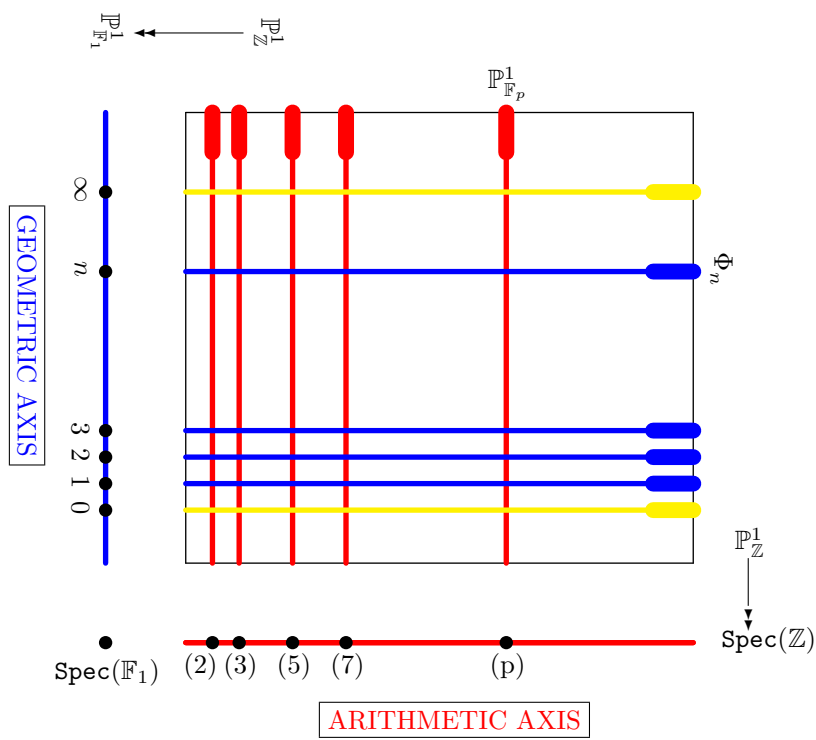

\[
\rotatebox{0}{
\begin{mfpic}[1]{0}{0}{100}{100}
\lines{(10,30),(200,30)}
\lines{(10,30),(10,200)}
\lines{(10,200),(200,200)}
\lines{(200,30),(200,200)}
\pen{2pt}
\drawcolor{red}
\lines{(10,0),(200,0)}
\lines{(20,30),(20,200)}
\lines{(30,30),(30,200)}
\lines{(50,30),(50,200)}
\lines{(70,30),(70,200)}
\lines{(130,30),(130,200)}
\lines{(10,0),(200,0)}
\pen{6pt}
\lines{(20,185),(20,200)}
\lines{(30,185),(30,200)}
\lines{(50,185),(50,200)}
\lines{(70,185),(70,200)}
\lines{(130,185),(130,200)}
\pen{2pt}
\drawcolor{yellow}
\lines{(10,50),(200,50)}
\lines{(10,170),(200,170)}
\pen{6pt}
\lines{(185,50),(200,50)}
\lines{(185,170),(200,170)}
\pen{2pt}
\drawcolor{blue}
\lines{(10,70),(200,70)}
\lines{(10,60),(200,60)}
\lines{(10,80),(200,80)}
\lines{(10,140),(200,140)}
\lines{(-20,30),(-20,200)}
\pen{6pt}
\lines{(185,70),(200,70)}
\lines{(185,60),(200,60)}
\lines{(185,80),(200,80)}
\lines{(185,140),(200,140)}
\pen{2pt}
\pointsize=5pt
\point{(20,0),(30,0),(50,0),(70,0),(130,0),(-20,50),(-20,60),(-20,70),(-20,80),(-20,140),(-20,170),(-20,0)}
\tlabel[cc](15,-10){(2)}
\tlabel[cc](30,-10){(3)}
\tlabel[cc](50,-10){(5)}
\tlabel[cc](70,-10){(7)}
\tlabel[cc](130,-10){(p)}
\tlabel[cc](-20,-12){$\wis{Spec}(\mathbb{F}_1)$}
\tlabel[cc](130,210){$\PP^1_{\F_p}$}
\tlabel[cl](215,37){$\PP^1_{\Z}$}
\tlabel[cc](220,20){$\begin{diagram} \\ \dOnto  \end{diagram}$}
\tlabel[cl](210,0){$\wis{Spec}(\Z)$}
\tlabel[cc](110,-30){\fbox{\textcolor{red}{ARITHMETIC AXIS}}}
\tlabel[cc](-50,115){\rotatebox{270}{\fbox{\textcolor{blue}{GEOMETRIC AXIS}}}}
\tlabel[cc](-30,50){\rotatebox{270}{$0$}}
\tlabel[cc](-30,60){\rotatebox{270}{$1$}}
\tlabel[cc](-30,70){\rotatebox{270}{$2$}}
\tlabel[cc](-30,80){\rotatebox{270}{$3$}}
\tlabel[cc](-30,140){\rotatebox{270}{$n$}}
\tlabel[cc](-30,170){\rotatebox{270}{$\infty$}}
\tlabel[cc](210,140){\rotatebox{270}{$\Phi_n$}}
\tlabel[cc](10,230){\rotatebox{270}{$\begin{diagram} \PP^1_{\Z}  \\ \dOnto \\ \PP^1_{\F_1} \end{diagram}$}}
\end{mfpic}}
\]
\caption{The extended Mumford drawing}
\end{figure}
But, you may wonder, what is this elusive 'field with one element' $\mathbb{F}_1$, and what do we mean with geometric objects defined over it? 

Two papers mark the beginning of this subject, one by Alexander L. Smirnov, "Hurwitz inequalities for number fields" \cite{Smirnov92} and one by Yuri I. Manin, "Lectures on zeta functions and motives" \cite{Manin92}, both originating from the fall of 1991. But, you will search in vein for any mention of the 'field with one element' in either paper. Manin introduces the 'absolute motive' and Smirnov 'an object that partially replaces the projective line over the constant field' in the number field case. Information about the historical origins and motivation behind these two papers can be gleaned from two letters from A.L. Smirnov to Y.I. Manin \cite{SmirnovLetters03} and an unpublished preprint by Smirnov and Mikhail Kapranov \cite{SmirnovKapranov96}.

It is an old idea to interpret the combinatorics of finite sets as the limit case of linear algebra over finite fields $\mathbb{F}_q$ when $q$ goes to $1$. This led to a folklore object, the 'field with one element', $\mathbb{F}_1$, vector spaces over which are just sets $V$ or pointed sets $V_{\bullet} = V \cup \{ 0 \}$ if we want to add a zero-vector. The dimension of the vector space $V$ is its cardinality $\# V$ and an $\mathbb{F}_1$-linear map between vector spaces $V \rTo W$ is just a map of sets (or a map of pointed sets $V_{\bullet} \rTo W_{\bullet}$ mapping the distinguished element of $V$ (the zero vector) to that of $W$). Consequently, one should interpret the general linear group $GL_n(\mathbb{F}_1)$ as the group of all permutations on a set with $n$ elements, that is $GL_n(\mathbb{F}_1) = S_n$. The analog of the usual determinant $det~:~GL_n(\mathbb{F}_1) \rTo \mathbb{F}_1^*$ is then the sign groupmorphism $sgn~:~S_n \rTo \{ \pm1 \}$ and hence one should view the alternating group $A_n$ as the special linear group $SL_n(\mathbb{F}_1)$. As an example of the slogan that linear algebra over $\mathbb{F}_1$ is the same as the combinatorics of finite sets consider an $n$-dimensional vector space $V/\mathbb{F}_1$ (that is, $\# V =n$), then the $k$-element subsets of $V$ should be viewed as points in the Grassmannian $\wis{Grass}(k,n)(\mathbb{F}_1)$ and hence its cardinality is the limit of the cardinalities of the actual Grassmannians $\wis{Grass}(n,k)(\mathbb{F}_q)$ over actual finite fields
\[
\begin{pmatrix} n \\ k \end{pmatrix} = \underset{q \rightarrow 1}{lim}~\#~\wis{Grass}(k,n)(\mathbb{F}_q) \]
One can play for some time exploring similar analogies, but quickly one feels that the setting lacks flexibility. In order to resolve this Alexander L. Smirnov introduced finite field extensions of $\mathbb{F}_1$. By analogy with the genuine finite field case one should have just one field-extension $\mathbb{F}_{1^n}$ of degree $n$ up to isomorphism and Smirnov proposed to take the monoid
\[
\mathbb{F}_{1^n} = \{ 0 \} \cup \pmb{\mu}_n \]
where $\pmb{\mu}_n$ is the group of all roots of unity of order $n$. 
Vector spaces over $\mathbb{F}_{1^n}$ can then be taken to be pointed sets $V_{\bullet} = \{ 0 \} \cup V$ where $V$ is a set having a free action of $\pmb{\mu}_n$, the dimension being the number of $\pmb{\mu}_n$-orbits. Linear maps are then maps of pointed sets which are also maps of $\mathbb{\mu}_n$-sets. Therefore, the Galois group $Gal(\mathbb{F}_{1^n}/\mathbb{F}_1)$ should be the multiplicative group $(\Z/n\Z)^*$ consisting of the power maps $\mu \mapsto \mu^d$ for all $(d,n)=1$. Taking limits, the algebraic closure of the field with one element should be considered as the pointed set,or monoid
$\overline{\mathbb{F}_1} = \{ 0 \} \cup \pmb{\mu}_{\infty}$
consisting of zero together with all roots of unity.

Smirnov and Kapranov remark in \cite{SmirnovKapranov96} that the idea of adjoining roots of unity as analogously to extension of the base field at least goes back to a letter by Andre Weil to Emil Artin \cite{Weil42} in which Weil writes : "ÒOur proof of the Riemann hypothesis (in the function field case) depended upon the extension of the function-fields by roots of unity, i.e. by constants; the way in which the Galois group of such extensions operates on the classes of divisors in the original field and its extensions gives a linear operator, the characteristic roots (i.e. the eigenvalues) of which are the roots of the zeta-function.
On a number field, the nearest we can get to this is by adjunction of $l$-th roots of unity, $l$ being fixed; the Galois group of this infinite extension is cyclic, and defines a linear operator on the projective limit of the (absolute) class groups of those successive finite extensions; this should have something to do with the roots of the zeta-function of the field. However, our extensions are ramified (but only at a finite number of places, viz. the prime divisors of $l$). Thus a preliminary study of similar problems in function-fields might enable one to guess what will happen in number-fields."

Smirnov's proposal was then to take as the schematic points of $\PP^1_{\F_1}$ to be the set
\[
\{ [0],[\infty ] \} \cup \{ [1],[2],[3],[4],[5],\hdots \}
\]
the set of all positive integers $\mathbb{N}$ together with a point at infinity. He also declares the {\em degree} of the point $n \in \mathbb{N}_0$ to be the Euler function $\phi(n)$, whereas the points $[0]$ and $[\infty ]$ have degree one. Here's why.
The {\em geometric} points of $\mathbb{P}^1_{\mathbb{F}_p}$ are of course
\[
\mathbb{P}^1_{\mathbb{F}_p}(\overline{\mathbb{F}}_p) = \{ [0]=[0:1],[\infty ]=[1:0] \} \cup \{ [\alpha ] = [\alpha : 1]~:~\alpha \in \overline{\mathbb{F}}^*_p \} \]
The Galois group $Gal(\overline{\mathbb{F}}_p/\mathbb{F}_p) = \hat{\mathbb{Z}}_+$ acts on this set by fixing the points $[0]$ and $[\infty ]$ and sending $\sigma([\alpha ] = [\sigma(\alpha)]$. The {\em schematic} points of $\mathbb{P}^1_{\mathbb{F}_p}$ are then the Galois-orbits for this action and the degree of a schematic point is the size of the corresponding orbit. If we assign to an orbit $\Oscr$ the polynomial
\[
\prod_{[\alpha ] \in \Oscr} (x-\alpha) \]
we see that the schematic points of $\mathbb{P}^1_{\mathbb{F}_p}$ consist of $[\infty]$ together with the set of all monic irreducible polynomials in $\mathbb{F}_p[x]$ and that the notion of degree of the schematic point coincides with the usual degree of the corresponding polynomial.

Here's an alternative description. We claim that we can identify the multiplicative group of the non-zero elements of the algebraic closure
\[
\overline{\mathbb{F}}^*_p \simeq \pmb{\mu}^{(p)} \]
with the group $\pmb{\mu}^{(p)}$ of all roots of unity having order prime to $p$. Clearly, if $\alpha \in \mathbb{F}_{p^n}^*$ then the  order of $\alpha$ is a divisor of $p^n-1$ and hence a number prime to $p$. Conversely, if $(m,p)=1$ then the residue class $\overline{p} \in \mathbb{Z}/m\Z$ is a unit and therefore for some integer $n$ we must have $p^n \equiv 1~\wis{mod}~m$. But then, $m | p^n -1$ and the primitive $m$-th roots of unity can be identified with a subgroup of the multiplicative group $\mathbb{F}^*_{p^n}$.
However, describing the correspondence $\overline{\mathbb{F}}^*_p \leftrightarrow \pmb{\mu}^{(p)}$ explicitly from a given construction $\overline{\mathbb{F}}_p$ is very challenging and we will address it later.

By {\em analogy} we can therefore define the {\em geometric} points of $\mathbb{P}^1_{\mathbb{F}_1}$ with
\[
\mathbb{P}^1_{\mathbb{F}_1}(\overline{\mathbb{F}}_1) = \{ [0],[\infty ] \} \cup \pmb{\mu}^{(1)}
\]
with $\pmb{\mu}^{(1)}$ the  set of all roots of unity with order prime to $1$, that is, the group $\pmb{\mu}_{\infty}$ of all roots of unity, leading to the proposal that 
\[
\overline{\mathbb{F}}_1 = \{ 0 \} \cup \pmb{\mu}_{\infty} \]
The schematic points of $\mathbb{P}^1_{\mathbb{F}_1}$ are then the  orbits of this set under the  action of the Galois group $Gal(\Q(\pmb{\mu}_{\infty})/\Q) = \hat{\Z}^*$. Clearly, these orbits are classified by
\[
\{ [0],[\infty ] \} \cup \{ [1],[2],[3],[4],[5],\hdots \} \]
where $[n]$ is the orbit consisting of all primitive $n$-th roots of unity, and hence the degree of the schematic point $[n]$ must be equal to the number of primitive $n$-th roots of unity, that is, to $\phi(n)$.

\section{The Habiro topology on $\PP^1_{\F_1}$}

The additive structure of the {\em profinite integers} \index{profinite integers} $\hat{\Z}$ as well as its multiplicative group of units $\hat{\Z}^*$ have already made their appearance as (absolute) Galois groups. As we will encounter them often, let us formally define these profinite integers following Hendrik Lenstra's account in \cite{Lenstra05}. 

Recall that any positive integer $n$ has a unique representation of the form
\[
n = c_k. k! + c_{k-1}.(k-1)! + \hdots + c_2.2! + c_1.1! \]
where the 'digits' $c_i$ are integers such that $0 \leq c_i \leq i$ for all $1 \leq i \leq k$ and $c_k \not= 0$. We then write $n$ in the {\em factorial number system} as $n = (c_kc_{k-1} \hdots c_2c_1)_!$, so for example $25=(1001)_!$. Profinite integers arise if we allow the sequences of digits to extend indefinitely to the left to get expressions such as $( \hdots c_4c_3c_2c_1)_!$. One can then identify positive integers $\N$ to be those profinite integers with $c_i=0$ for $i>>$. Also the negative integers $-\N_0$ can be characterized as those profinite integers such that $c_i=i$ for all but finitely many $i$. For example, $-1 = (\hdots 654321)_!$, that is $c_i=i$ for all $i$. 

To add two profinite integers, one adds them digitwise, proceeding from the right and if the sum of the two $i$-th digits is larger than $i$, one substracts $i+1$ from it and adds a carry of $1$ to the sum of the $i+1$-th digits. With this rule we have indeed that $-1$ (as above) $+1=0$. To multiply two profinite integers we use the rule that the first $k$ digits of the product $s \times t$ depend only on the first $k$ digits of $s$ and $t$, thereby reducing the problem of computing products to the case of ordinary positive integers. These operations make the profinite integers $\hat{\Z}$ into a commutative ring with unit element $1$. Those in the know will have observed that all we did was work out the ring-rules for the projective limit $\hat{\Z} = \underset{\leftarrow}{lim}~\Z/n! \Z$. 

But let us return to $\PP^1_{\F_1}$. We have defined, following Smirnov, that the schematic points of $\PP^1_{\F_1}$ are the orbits of $\pmb{\mu}_{\infty}$ under the action of the multiplicative group of profinite integers $\hat{\Z}^*$ which is the Abelian Galois group $Gal(\Q(\pmb{\mu}_{\infty})/\Q)$. We would now like to define a topology on $\pmb{\mu}_{\infty}$ compatible with this action, and clearly the induced analytic topology does not satisfy this condition.

In \cite{Habiro08} Kazuo Habiro introduced a new topology on $\pmb{\mu}_{\infty}$ in order to unify invariants of $3$-dimensional homology spheres, introduced first by Ed Witten by means of path integrals and rigorously constructed by Reshitikhin and Turaev. Two roots of unity $\alpha,\beta \in \pmb{\mu}_{\infty}$ are said to be {\em adjacent} if their quotient $\alpha \beta^{-1}$ is of pure prime power order $p^m$ for $m \in \Z$ and $p$ a prime number, or equivalently, if the difference $\alpha-\beta$ is not a unit in the integral closure of $\Z$ in $\Q(\alpha,\beta)$. Clearly, the action of the absolute Galois group $Gal(\overline{\Q}/\Q)$ and of $\hat{\Z}^* = Gal(\Q(\pmb{\mu}_{\infty})/\Q)$ preserves adjacency. 
\[
\includegraphics[height=4.1cm]{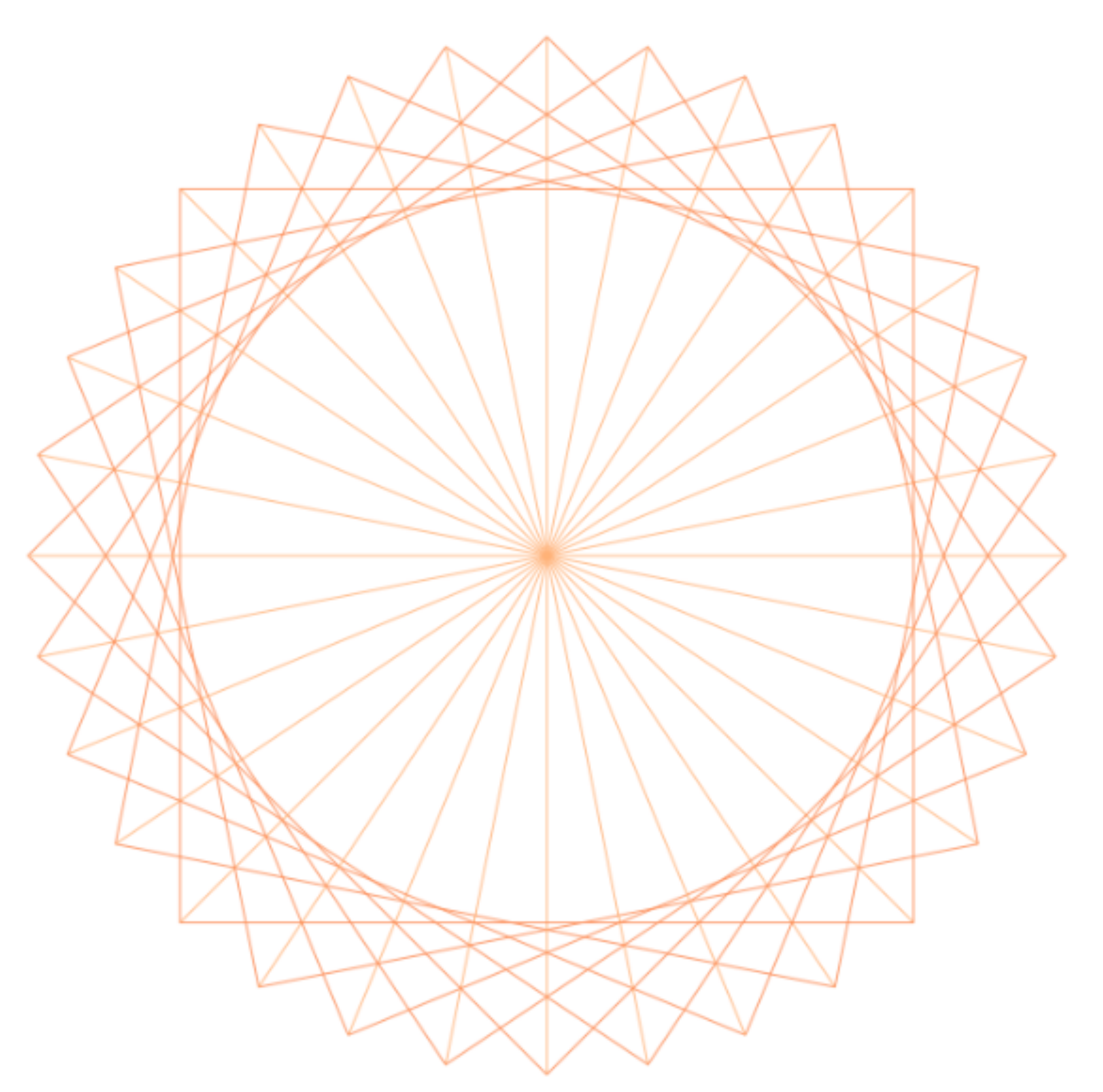}
\includegraphics[height=4.1cm]{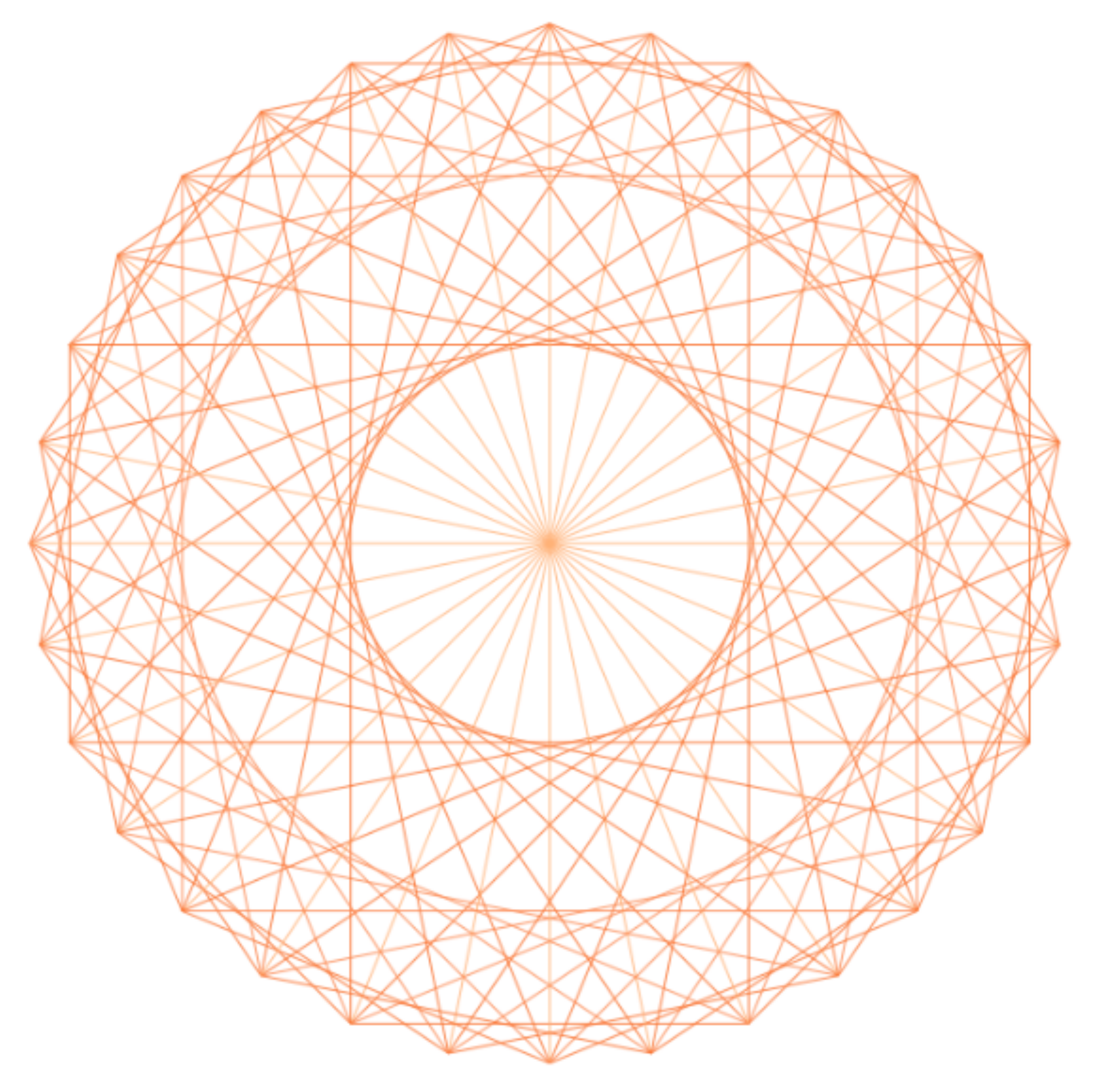}
\includegraphics[height=4.1cm]{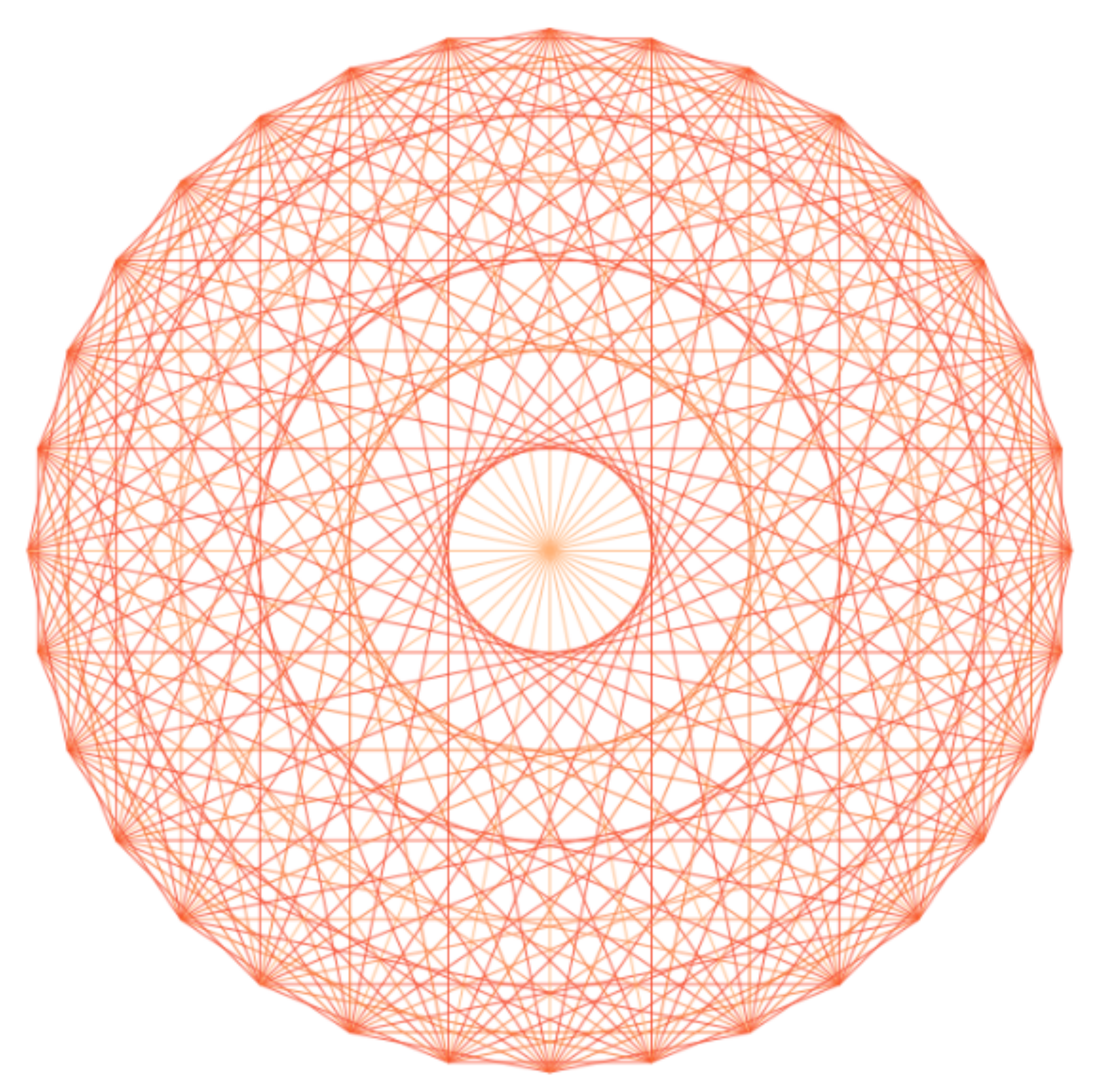}
\]
\[
\includegraphics[height=7cm]{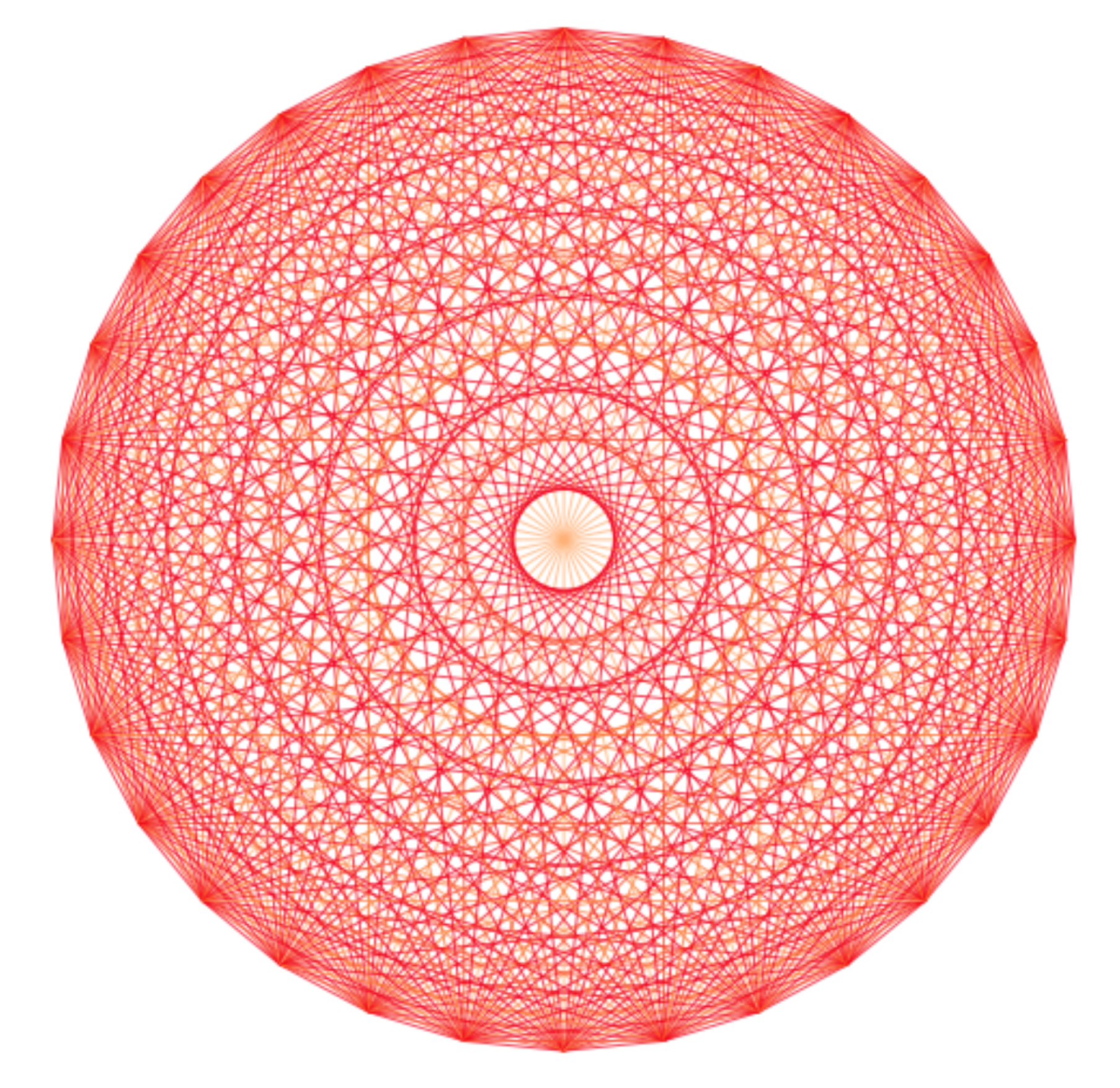} \]
Above we indicated adjacency of $32$-nd roots of unity under increasing powers $2^k$, respectively $k \leq 2,3,4$ and finally $5$.

The {\em Habiro topology} \index{Habiro topology} on $\pmb{\mu}_{\infty}$ is then defined by taking as its open sets those subsets $U \subset \pmb{\mu}_{\infty}$ having the property that for every $\alpha \in U$ all except finitely many $\beta \in \pmb{\mu}_{\infty}$ adjacent to $\alpha$ also belong to $U$. The Galois action is continuous in this topology which is in sharp contrast to the induced analytic topology. Further, any cofinite subset is clearly open in the Habiro topology, but we will soon see that there are plenty of other open subsets. In fact we will show that the Habiro topology is not locally compact on $\pmb{\mu}_{\infty}$. 

But let us first explain the $\F_1$-origins of the Habiro topology. Our depiction of $\PP^1_{\Z}$ as the product of the arithmetic and the geometric axis was an over-simplification. Whereas the vertical fibers $\mathbb{V}(p)$ (the red lines) are clearly disjoint, this is not necessarily the case for the blue lines $\mathbb{V}(\Phi_n)$, because the height one prime ideals $(\Phi_n(x))$ and $(\Phi_m(x))$ do not have to be comaximal and whenever $(\Phi_n(x),\Phi_m(x)) \not= \Z[x]$ there will be a point lying on both curves, so the 'lines' will intersect over certain prime numbers. This must happen as we know that the prime spectrum of the integral group ring $\wis{Spec} \Z C_n$ is connected, see for example \cite{Serre77}. Hence, its minimal prime ideals which are all of the form $\mathbb{V}(\Phi_d(x))$ for $d | n$ will intersect. For cyclotomic polynomials we have complete information about potential comaximality:
\begin{itemize}
\item{If $\tfrac{m}{n} \not= p^k$ for some prime number $p$, then $(\Phi_m(x),\Phi_n(x))=\Z[x]$ so these cyclotomic prime ideals are comaximal and the corresponding blue lines will not intersect, however}
\item{if $\tfrac{m}{n} = p^k$ for some prime number $p$, then $\Phi_m(x) \equiv \Phi_n(x)^d~\wis{mod}~p$ for some integer $d$ and hence these cyclotomic prime ideals are not comaximal and the corresponding blue lines will intersect over $p$.}
\end{itemize}
In the second case we have non-split extensions as $\Z[x,x^{-1}]$-modules
\[
\begin{cases}
0 \rTo \Z[\zeta_n] \rTo E \rTo \Z[\zeta_m] \rTo 0 \\
0 \rTo \Z[\zeta_m] \rTo F \rTo \Z[\zeta_n] \rTo 0 
\end{cases}
\]
In fact, Fritz-Erdmann Diederichsen, a student of Zassenhaus, calculated already in 1940 that there are exactly $p^{\mu(min(m,n))}$ such extensions, in either direction \cite{Diederichsen40} (see also \cite[Thm. 25.26]{CurtisReiner81}). 

In noncommutative algebraic geometry we are very familiar with such situations. If $S$ and $T$ are two finite dimensional simple representations of a $\C$-algebra $A$ such that $Ext^1_A(S,T) \not= 0$, we say that their annihilating maximal ideals belong to the same 'clique' and know that we should think of $S$ and $T$ as two noncommutative points lying infinitely close together, or equivalently, that $S$ and $T$ share some tangent information.
Using this noncommutative intuition we therefore define a clique- or adjacency relation on pairs of natural numbers
\[
m \sim n \qquad \text{if and only if} \qquad \frac{m}{n} = p^{\pm a} \]
for some prime number $p$. It is this clique-relation which lies behind the definition of the Habiro topology on $\pmb{\mu}_{\infty}$.
Below, we depict the inter-weaving patters of the horizontal lines $\Phi_1, \Phi_p, \Phi_q, \Phi_{p^2}$ and $\Phi_{pq}$ for prime numbers $p < q$.
\[
\begin{mfpic}[1]{0}{0}{100}{100}
\lines{(0,0),(200,0)}
\lines{(0,0),(0,200)}
\lines{(0,200),(200,200)}
\lines{(200,0),(200,200)}
\lines{(50,0),(50,200)}
\lines{(80,0),(80,200)}
\pen{2pt}
\drawcolor{yellow}
\lines{(0,0),(200,200)}
\drawcolor{red}
\lines{(0,50),(50,50)}
\lines{(50,50),(200,100)}
\drawcolor{blue}
\lines{(0,80),(130,80)}
\lines{(150,80),(200,80)}
\drawcolor{green}
\lines{(0,100),(15,85)}
\lines{(25,75),(50,50)}
\lines{(50,50),(200,50)}
\drawcolor{black}
\lines{(0,120),(50,80)}
\lines{(50,80),(62,72)}
\lines{(71,66),(80,60)}
\lines{(80,60),(122,53)}
\lines{(158,47),(200,40)}
\tlabel[cc](-10,50){$p$}
\tlabel[cc](-10,80){$q$}
\tlabel[cc](-10,100){$p^2$}
\tlabel[cc](-10,120){$pq$}
\tlabel[cc](-10,0){$1$}
\tlabel[cc](50,-10){$(p)$}
\tlabel[cc](80,-10){$(q)$}
\end{mfpic}
\]
The {\em Habiro topology} \index{Habiro topology} on $\PP^1_{\F_1}= \{ [n]~:~n \in \N_0^{\times} \} \cup \{ [0],[\infty ] \}$ is therefore defined as follows.  An open set is then a subset of $\PP^1_{\F_1}$ of the form
\[
U \quad \text{or} \quad U^0=U \cup \{ [0] \} \quad \text{or} \quad U^{\infty}=U \cup \{ [ \infty ] \} \quad \text{or} \quad U^{0\infty}=U \cup \{ [0],[ \infty ] \} \]
where $U$ has the property that if $[m] \in U$ then all but finitely $[n]$ such that $m \sim n$ also belong to $U$. Clearly, all cofinite subsets of $\PP^1_{\F_1}$ are open, but there are more. For a prime number $p$ define the set
\[
U_p = \{ [n]~|~\text{if $p | n$ then $p^2 | n$} \} \]
then $U_p$ is open for if $[n] \in U_p$ and $p \nmid n$ then there is all $[m] \in U_p$ if $n \sim m$ except for the one point $[ pn ]$ and if $n = p^k n'$ then again all $[m] \in U_p$ when $m \sim n$ except for the one point $[ p n' ]$. But still, the complement of $U_p$ is infinite as
\[
\PP^1_{\F_1} - U_p = \{ [ pa ]~|~(p,a)=1 \} \cup \{ [0],[\infty ] \} \]
More generally, if $m = p_1^{k_1} \hdots p_l^{k_l}$ then the set
\[
U_m = \{ [n]~|~\forall~1 \leq i \leq l~:~\text{if $p_i | n$ then $p_i^{k_i+1} | n$} \} \]
is open in the Habiro topology and we have relations such as $U_m \cap U_n = U_{lcm(m,n)}$. 

A striking feature is that $\PP^1_{\F_1}$ is {\em not} compact in the Habiro topology. Indeed, as $[n] \in U_q$ for all primes $q$ such that $q \nmid n$ we have that
\[
\PP^1_{\F_1} = \bigcup_{p~\text{prime}}~U_p^{0\infty} \]
but no finite sub-cover exists as $[ p_1\hdots p_k ] \notin U_{p_1}^{0 \infty} \cup \hdots \cup U_{p_k}^{0 \infty}$.

\section{The \'etale site of $\F_1$}

Until now our exposition has been pretty intuitive and we would like to have a solid framework to give formal definitions for these elusive objects as well as to perform actual calculations. Such a proposal was put forward by Jim Borger in 2009 in \cite{Borger09} and streamlined a plethora of previous attempts, all giving more or less the same class of examples. For a map of $\F_1$-land before Borger we refer to \cite{PenaLorscheid09}.

Borger's approach is based on the notion of {\em $\lambda$-rings} \index{$\lambda$-rings} in the sense of Grothendieck's Riemann-Roch theory \cite{Grothendieck58}. For us, a {\em $\lambda$-ring} will be a commutative ring $R$ with identity {\em without additive torsion} equipped with a collection of maps $\lambda^n~:~R \rTo R$ for $n \in \N$ satisfying an unwieldy list of axioms, see for example \cite[Chapter 1]{Knutson73}. By \cite{Wilkerson82}, assuming our no-torsion assumption, $R$ is a $\lambda$-ring if and only if there is a collection of {\em commuting} ring-endomorphisms
$\psi^p~:~R \rTo R$ for all prime numbers $p$
such that $\psi^p$ is a lift of the Frobenius-automorphism on $R \otimes_{\Z} \mathbb{F}_p$, that is, we have a commuting diagram for all prime numbers $p$
\[
\begin{diagram}
R & \rTo^{\psi^p} & R \\
\dOnto & & \dOnto \\
R/pR & \rTo^{Frob} & R/pR
\end{diagram}
\]
Borger's proposal is to {\em define} $\mathbb{F}_1$-algebras as $\Z$-algebras with a $\lambda$-ring structure, the idea being that one can interpret the collection of commuting endomorphisms $\{ \psi^p \}$ as descent data from $\Z$-algebras to $\mathbb{F}_1$-algebras. $\mathbb{F}_1$-algebra maps are then morphisms of $\lambda$-rings, that is, ringmorphisms commuting with the Frobenius lifts $\psi^p$. Accepting this proposal it then follows that base ring extension $- \otimes_{\mathbb{F}_1} \Z$ can be viewed formally as the operation of forgetting the $\lambda$-ring structure, that is, of stripping off the $\F_1$-structure.

In this way we can {\em define} $\mathbb{F}_1$ to be $\Z$ with its unique $\lambda$-ring structure in which all $\psi^p = id_{\Z}$. Similarly, the polynomial ring $\mathbb{F}_1[x]$ can be defined to be the integral polynomial ring $\Z[x]$ equipped with the $\lambda$-ring structure defined by $\psi^p(x) = x^p$ which is indeed a $\Z$-endomorphism lift of the Frobenius automorphism on $\mathbb{F}_p[x]$ by little Fermat and the binomial theorem. The field $\mathbb{F}_{1^n}$ should then be taken to be the integral groupring $\Z[\pmb{\mu}_n]$ with  $\lambda$-rings structure induced by $\psi^p(\mu) = \mu^p$. This gives the fanciful identities
\[
\mathbb{F}_{1^n} \otimes_{\mathbb{F}_1} \Z = \Z[\pmb{\mu}_n] \qquad \mathbb{F}_1[x] \otimes_{\mathbb{F}_1} \Z = \Z[x] \quad \text{and} \quad \PP^1_{\F_1} \otimes_{\F_1} \Z = \PP^1_{\Z} \]
where on the right hand sides we forget about the $\lambda$-structures.

If $G$ is a finite group and let $\chi_1,\hdots,\chi_h$ be its irreducible characters, then the Grothendieck ring of finite dimensional $\C[G]$-representations
\[
R(G) = \Z \chi_1 \oplus \hdots \oplus \Z \chi_h \]
has a $\lambda$-ring structure defined by $\psi^p(\chi) = \chi'$ where $\chi'$ is the class function $\chi'(g) = \chi(g^p)$ for all $g \in G$, see for example \cite[\S 9.1]{Serre77}. If $p$ does not divide the order of $G$, then $\psi^p$ is an automorphism permuting the irreducible characters. As the $\psi^p$ commute we can define operations $\psi^n = (\psi^{p_1})^{\circ a_1} \circ \hdots \circ (\psi^{p_k})^{\circ a_k}$ if $n=p_1^{a_1} \hdots p_k^{n_k}$. The $\psi^n$ are the Adams operations and they give an action of the multiplicative monoid $\N^{\times}_0$ on $R(G)$. On the other hand we have that
\[
\Q(\pmb{\mu}_{\infty})\otimes_{\Z} R(G) \simeq \underbrace{\Q(\pmb{\mu}_{\infty}) \times \hdots \times \Q(\pmb{\mu}_{\infty})}_h \]
and hence there is a Galois action by $\hat{\Z}^* = Gal(\Q(\pmb{\mu}_{\infty})/\Q)$ on $\Q(\pmb{\mu}_{\infty})\otimes_{\Z} R(G)$.
This Galois action is compatible with the action by Adams operations on $R(G)$ via the embedding $\N^{\times}_0 \rInto \hat{\Z}^*$, see for example \cite[\S 12.4]{Serre77}.

James Borger and Bart De Smit vastly generalized this example in \cite{BorgerDeSmit08} to include all reduced $\lambda$-rings $R$ which are finite projective $\Z$-modules, allowing us to describe the Galois (or \'etale) site of $\mathbb{F}_1$. 

Let us first consider the case when $K = R \otimes_{\Z} \Q$ is a $\lambda$-ring and a finite \'etale $\Q$-algebra. By Grothendieck's version of Galois theory we have an anti-equivalence of categories between the category of finite \'etale $\Q$-algebras and the category of finite discrete sets equipped with a continuous action of the absolute Galois group $G=Gal(\overline{\Q}/\Q)$, assigning to $K$ the set
\[
S = \wis{alg}_{\Q}(K,\overline{\Q}) \]
on which $\sigma \in G$ acts by left composition, that is, $\sigma.s = \sigma \circ s$. If we have in addition a $\lambda$-ring structure, that is a commuting family of $\Q$-endomorphisms $\Psi^n$ on $K$, then we have also an action of the monoid $\N^{\times}_0$ on $S$ acting by composition on the right, that is $n.s = s \circ \Psi^n$. Consequently, the category of rational $\lambda$-algebras which are finite \'etale over $\Q$ is the category of finite discrete sets equipped with a continuous action of the monoid $Gal(\overline{\Q}/\Q) \times \N^{\times}_0$ where $\N^{\times}_0$ is given the discrete topology. 

The Galois action on $S$ gives us a groupmorphism
\[
Gal(\overline{\Q/}\Q) \rTo \wis{perm}(S) \subset \wis{maps}(S) \]
where $\wis{perm}(S)$ are all permutations on $S$ and $\wis{maps}(S)$ are all set-maps from $S$ to itself. The kernel $N$ gives us a finite Galois extension $L = \overline{\Q}^N$ of $\Q$ with Galois group $\overline{G} = Gal(\overline{\Q}/\Q)/N$ and as $R \otimes_{\Z} \Q = K = L_1 \times \hdots \times L_k$, all $L_i$ are subfields of $L$ and hence every $\sigma \in \overline{G}$ induces an automorphism on $L_i$. Let $\Oscr_L$ be the integral closure of $\Z$ in $L$, then by Galois theory we have
\[
S = \wis{alg}_{\Q}(K,\overline{\Q}) = \wis{alg}_{\Q}(K,L) = \wis{alg}_{\Z}(R,\Oscr_L) \]
and we have the action by $(\sigma,n) \in \overline{G} \times \N_0^{\times}$ on $S$ via
\[
\xymatrix{ R \ar[rr]^s \ar@{.>}[drr] & & \Oscr_L \ar[d]^{\sigma} \\
R \ar[u]^{\psi^n} \ar@{.>}[urr] & & \Oscr_L} \]
By Chebotarev's theorem there are infinitely many prime numbers $p$ such that there exists a prime ideal $P \triangleleft \Oscr_L$ lying over $p$ with $\sigma$ a lift of the Frobenius-automorphism on $\Oscr_L/P$, so we can choose a prime $p$ not dividing the discriminant $\Delta(R)$. Note that if $p \nmid \Delta(R)$ there is a unique lift of the Frobenius-map which is automatically an automorphism by the category equivalence between \'etale $\F_p$- and \'etale $\hat{\Z}_p$-algebras. But then, the restriction of $\sigma$ to $R$ via the embedding $s$ is equal to $\Psi^p$ and we have $\sigma \circ s = s \circ \Psi^p$ on $R$. As this holds for any $\sigma \in \overline{G}$ we have that the image of $\overline{G} \rTo \wis{perm}(S) \subset \wis{maps}(S)$ is contained in the image of $\N_0^{\times} \rTo \wis{maps}(S)$. But, as $\N_0^{\times}$ is an Abelian monoid, it follows that the image in $\wis{perm}(S)$ is Abelian and hence that $L/\Q$ is an Abelian Galois extension!

But now we can invoke the Kronecker-Weber theorem asserting that $L \subset \Q(\pmb{\mu}_c)$ where $c$ in only divisible by primes $p$ which ramify in $L$ and as $L$ is the common Galois extension of the components of $R \otimes_{\Z} \Q$ those $p$ must also divide $\Delta(R)$. That is, there exists $c \in \N$ with all its prime factors dividing the discriminant $\Delta(R)$ such that the Galois action on $S = \wis{alg}_{\Z}(R,\Oscr_L)$ factorizes through the cyclotomic character
\[
Gal(\overline{\Q}/\Q) \rOnto (\Z/c\Z)^* = Gal(\Q(\pmb{\mu}_c)/\Q) \]
and for all $p \nmid \Delta(R)$ the action op $p \in \N_0^{\times}$ on $S$ is equal to that of $p~\wis{mod}~c \in Gal(\Q(\pmb{\mu}_c)/\Q)$. 

So far, we have factorized the $Gal(\overline{\Q}/\Q) \times \N_0^{\times}$-action on $S$ via $\hat{\Z}^* \times \N_0^{\times}$ and now we want to factorize it further through the map $\hat{\Z}^* \times \N_0^{\times} \rTo \hat{\Z}^{\times}$ via the natural embeddings of both factors and where $\hat{\Z}^{\times}$ is the multiplicative monoid of the profinite integers $\hat{\Z}$.

We can apply the foregoing for every $d \in \N_0^{\times}$ on the sub $\lambda$-ring $\Psi^d(R) \subset R$ with corresponding finite set $d.S$. That is, there exists $c_d \in \N_0$ such that the Galois-action on $d.S$ factorizes through $(\Z/c_d\Z)^*$ and such that the action of any $n$ with $nd.S=d.S$ is the same as the action of $n~\wis{mod}~c_d$. For every prime number $p$ let $a_p$ be the smallest power such that $p^{a_p+1}.S = a^{a_p}.S$ and as $a_p > 0$ only for those $p | \Delta(R)$ we have a finite number $r_0=\prod_{p|\Delta(R)} p^{a_p}$ which satisfies the property that for all $n \in \N_0$ we have that $n.S = gcd(n,r_0).S$. Now let $r$ be the least common multiple of all $d.c_d$ where $d | r_0$, then we claim that the above action factorizes through the multiplicative monoid $(\Z/c_d \Z)^{\times}$, that is, we have to show that if $d_1~\equiv~d_2~\wis{mod}~r$ that the actions of $d_1$ and $d_2$ on $S$ coincide. As $r_0 | r$ we have $gcd(d_1,r_0)=gcd(d_2,r_0)=d$ whence $d_1.S=d.S=d_2.S$. If we write $d_i=d d'_i$ this entails that $(d_i',c_d)=1$ but then $d_1=dd'_1 \equiv dd'_2=d_2~\wis{mod}~dc_d$. But then, $d'_1 \equiv d'_2~\wis{mod}~c_d$ and so they act the same on $d.S$ whence $d_1$ and $d_2$ act the same on $S$!

This then is the main result of Borger and De Smit \cite{BorgerDeSmit08} that a necessary and sufficient condition for the existence of an integral $\lambda$-ring $R$ which is finite and projective over $\Z$ contained in the $\lambda$-ring $R \otimes_{\Z} \Q = K$ is that the action of the monoid $Gal(\overline{\Q}/\Q) \times \N^{\times}_0$ on the finite set $S$ describing $K$ factors through an action of the monoid $\hat{\Z}^{\times}$. It follows that the category of all such $\lambda$-rings is anti-equivalent to the category of finite discrete sets with a continuous action of the monoid $\hat{\Z}^{\times}$ and that every such $\lambda$-ring is contained as $\lambda$-ring in a product of cyclotomic fields, generalizing the case of the Grothendieck ring $R(G)$ of a finite group $G$. 

Motivated by Grothendieck's interpretation of Galois theory we have the fanciful picture of the {\em absolute Galois monoid} \index{absolute Galois monoid} of the field with one element $\mathbb{F}_1$
\[
Gal(\overline{\mathbb{F}_1}/\mathbb{F}_1) \simeq \hat{\Z}^{\times} \]
Because the subset $0.S \subset S$ is Galois invariant it corresponds to a factor $\Q$ in $K$, so $K$ can never be a field unless $K=\Q$. In particular, whereas $\Z[\pmb{\mu}_n]$ is an integral $\lambda$-ring, the subring $\Z[\zeta_n]$ of the cyclotomic field $\Q(\zeta_n)$ where $\zeta_n$ is a primitive $n$-th root of unity is {\em not}.

For example, let us work out the $\lambda$-ring structure on $R(S_3)$, the representation ring of the symmetric group $S_3$ and its associated finite $\hat{\Z}$-set. For any finite group $G$ let $X$ be the set of conjugacy classes in $G$, then we can identify this set with
\[
S = \wis{alg}_{\Z}(R(G),\C) = \{ x~:~R(G) \rTo \C~|~x(V_i) = \chi_{V_i}(x)~\forall~V_i \in \wis{irreps}(G) \} \]
Moreover, one knows in general that the discriminant $\Delta(R(G)) = \frac{(\# G)^{\# X}}{\prod_{x \in X} \# x}$. Specializing to the case when $G=S_3$ we have that $\Delta(R(S_3)) = 36$ and recall that the character table of $S_3$ is
\[
\begin{array}{c|ccc}
x= & [1] & [2] & [3] \\
& () & (1,2) & (1,2,3) \\
\hline 
V_1 & 1 & 1 & 1 \\
V_2 & 1 & -1 & 1 \\
V_3 & 2 & 0 & -1
\end{array}
\]
The Frobenius lifts (aka Adams operators) send a class-function $\chi$ to the class-function $\Psi^n(\chi)(g)=\chi(g^n)$. As $()^n=()$ for all $n$, $(1,2)^n=()$ for even $n$ and $=(1,2)$ for odd $n$ and $(1,2,3)^n=()$ for $n$ a multiple of $3$ and is conjugated to $(1,2,3)$ otherwise. Therefore, if $\chi_i$ is the character-function of $V_i$ and computes from the character-table that for prime numbers $p$ we have
\begin{itemize}
\item{$\Psi^p(\chi_1) = \chi_1,~\forall p$}
\item{$\Psi^2(\chi_2) = \chi_1$ and $\Psi^p(\chi)=\chi_2,~\forall p \not=2$}
\item{$\Psi^2(\chi_3) = \chi_1+\chi_3-\chi_2$, $\Psi^3(\chi_3) = \chi_1+\chi_2$ and $\Psi^p(\chi_3) = \chi_3, \forall p \not= 2,3$}
\end{itemize}
which determines the $\lambda$-ring structure on $R(S_3)$. The action of $n \in \N_0^{\times}$ on the algebra map $[i] \in \wis{alg}_{\Z}(R(S_3),\C)$ is given by $n.[i] = [i] \circ \Psi^n$ and hence it follows from the above that $p.[i]=[i]$ for all primes $p \not= 2,3$ and one verifies that the action of $2$ and $3$ is given by the following maps on $S=\{ [1],[2],[3] \}$
\[
\xymatrix{& [1] \ar@[blue]@(r,ur)_3 \ar@[red]@(l,ul)^2 & \\
[2] \ar@[blue]@(dl,dr)_3 \ar@[red][ru]^2 & & [3] \ar@[red]@(dl,dr)_2 \ar@[blue][ul]_3 } \]
From this it follows that $2.S = 2^2.S = \{ [1],[3] \}$ and $3.S=3^2.S=\{ [1],[2] \}$ whereas $6.S=12.S=18.S=\{ [1] \}$. Further, the Galois action on $R(G)$ and any of its sub $\lambda$-rings is trivial. With the notations used before we therefore get  that $r_0=6$ and all $c_d=1$ showing that the $\hat{\Z}^{\times}$-action on $S$ factorizes through the multiplicative monoid action of $(\Z/6\Z)^{\times}$ as indicated in the above colored graph.

\section{What is $\PP^1_{\F_1}$?}

Now that we have a formal definition of $\F_1$-algebras, namely those $\Z$-rings without additive torsion which are $\lambda$-rings, it makes sense to define for any such $\lambda$-ring $R$ its {\em $\lambda$-spectrum} \index{$\lambda$-spectrum} which is the collection of all kernels of $\lambda$-ring morphisms from $R$ to reduced $\lambda$-rings
\[
\wis{Spec}_{\lambda}(R) = \{ \wis{ker}(R \rTo^{\phi} A)~|~\text{$A$ a reduced $\lambda$-ring and $\phi \in \wis{alg}_{\lambda}(R,A)$}~\} \]
which is clearly functorial. 

The geometric or $\overline{\F_1}$-points in the $\lambda$-spectrum then correspond to kernels of $\lambda$-ringmorphisms $R \rTo A$ where $A$ is one of the integral $\lambda$-rings described in the previous section, that is a finite projective $\Z$-ring with $\lambda$-structure such that $A \otimes_{\Z} \Q$ is an \'etale $\Q$-algebra. We have seen that such rings are of the form $A=A_S$ where $S$ is a finite set with a continuous monoid action by $\hat{\Z}^{\times}$. As these sets are ordered by inclusions $S \subset T$ compatible with the $\hat{\Z}^{\times}$-action and via the anti-equivalence this corresponds to $\lambda$-ring epimorphisms $A_T \rOnto A_S$, it makes sense to define the {\em maximal $\lambda$-spectrum} \index{maximal $\lambda$-spectrum} of $R$ to be
\[
\wis{max}_{\lambda}(R) = \{ \wis{ker}(R \rOnto^{\phi} A)~|~\text{$A$ \'etale over $\F_1$ and $\phi \in \wis{alg}_{\lambda}(R,A)$}~\} \]
As this space may still be too hard to compute in specific examples, we often reduce to the subset of all  {\em cyclotomic points} \index{cyclotomic points} (or in Manin-parlance of \cite{Manin08}, the {\em cyclotomic coordinates} \index{cyclotomic coordinates}) which is the set
\[
\wis{max}_{cycl}(R) = \{ \wis{ker}(R \rOnto^{\phi} \Z[\pmb{\mu}_n])~|~n \in \N,~\phi \in \wis{alg}_{\lambda}(R,\Z[\pmb{\mu}_n])~\} \]
For the $\lambda$-ring $\Z[x]$ we have that
\[
\wis{Spec}_{\lambda}(\Z[x]) = \{ 0 \} \cup \wis{max}_{\lambda}(\Z[x]) \quad \text{and} \quad \wis{max}_{cycl}(\Z[x]) = \{ (\Phi_n(x))~|~n \in \N \} \]
as any $\lambda$-ring epimorphism $\Z \rOnto \Z[\pmb{\mu}_n]$ must map $x$ to $x^i$ with $(i,n)=1$, that is, to a primitive $n$-th root of unity. Hence, we now finally have a formal definition of $\PP^1_{\F_1}$: it is the set of cyclotomic points of $\PP^1_{\Z}$, equipped with the Habiro topology.

One can again use methods from noncommutative algebraic geometry to obtain 'geometric objects' and their associated 'rings of functions' and apply this to the setting of $\mathbb{F}_1$-geometry to arrive at a similar description. 

In \cite{KontSoib06} Maxim Kontsevich and Yan Soibelman introduce a {\em noncommutative thin scheme} \index{noncommutative thin scheme} (over the complex numbers) to a covariant functor $\mathbb{X}~:~\wis{fd-alg}_{\C} \rTo \wis{sets}$ from finite dimensional (not necessarily commutative) $\C$-algebras to sets, commuting with finite projective limits. They show that such thin schemes are represented by a coalgebra $C_{\mathbb{X}}$ which they call the coalgebra of distributions on $\mathbb{X}$ and its dual algebra $C_{\mathbb{X}}^*$ is then called the algebra of function $\mathcal{O}(\mathbb{X})$ on $\mathbb{X}$. We will be interested in affine thin schemes, that is to a $\C$-algebra we associate its {\em representation functor}
\[
\wis{rep}_A~:~\wis{fd-alg}_{\C} \rTo \wis{sets} \qquad B \mapsto \wis{alg}_{\C}(A,C) \]
By Kostant duality (see for example \cite[Chapter VI]{Sweedler69}) this thin scheme is represented by the {\em dual coalgebra} \index{dual coalgebra} $A^o$ which consists of all linear functionals on $A$ which factor through a finite dimensional algebra quotient of $A$
\[
A^o = \{ f \in A^*~|~Ker(f) \supset I \triangleleft A~:~dim_{\C}(A/I) < \infty \} \]
and hence its corresponding ring of functions is $(A^o)^*$. One can use the $A_{\infty}$-structure on Yoneda-Ext algebras to describe the structure of the dual coalgebra $A^o$ for general $A$, see \cite{LeBruyn08}.

The motivating example being $X$ a commutative (complex) affine variety, when the dual coalgebra $\C[X]^o$ decomposes over the points of $X$ as distinct maximal ideals $\mathfrak{m}_x$ are comaximal
\[
C_X = \C[X]^o = \bigoplus_{x \in X} C_{X,x} \]
where $C_{X,x}$ is a subcoalgebra of the enveloping coalgebra $U(T_{X,x})$ of the Abelian Lie algebra on the Zariski tangent space $T_{X,x} = (\mathfrak{m}_x/\mathfrak{m}_x^2)^*$. Consequently, the ring of functions also decomposes over the points
\[
(\C[X]^o)^* = \prod_{x \in X} \hat{\Oscr}_{\mathfrak{m}} \]
where $\hat{\Oscr}_{\mathfrak{m}}$ is the $\mathfrak{m}_x$-adic completion of the local ring $\mathcal{O}_{\mathfrak{m}_x}$. Hence, the dual coalgebra contains a lot of geometric information: the points of $X$ can be recovered from it as the simple factors of the coradical $corad(\C[X]^o)$ and its dual algebra gives us the basics of the \'etale topology on $X$.  

Let us illustrate this in the case of interest, that is when $X = \mathbb{A}^1_{\C}$ with coordinate ring $\C[x]$. Every cofinite dimensional ideal of $\C[x]$ is of the form $I=((x-\alpha_1)^{e_1} \hdots (x-\alpha_k)^{e_k}))$ and as the different factors are comaximal, linear functionals on $\C[x]/I$ split over the distinct factors
\[
(\frac{\C[x]}{I})^* = (\frac{\C[x]}{((x-\alpha_1)^{e_1})})^* \oplus \hdots \oplus (\frac{\C[x]}{((x-\alpha_k)^{e_k})})^* \]
Each of these factors is the dual coalgebra of a truncated polynomial ring and if we take $z^i$ to be the basis dual to the $y^i$ we have
\[
(\frac{\C[y]}{(y^n)})^* = \C 1 + \C z + \hdots + \C z^{n-1} \quad \text{with} \quad \begin{cases} \Delta(z^k) = \sum_{i+j=k} z^i \otimes z^j \\ \epsilon(z^i) = \delta_{0i} \end{cases} \]
that is the structure of the truncated enveloping algebra. Hence we have proved that 
\[
\C[x]^o = \bigoplus_{\alpha \in \mathbb{A}^1_{\C}} U(T_{\mathbb{A}^1_{\C},\alpha}) \quad \text{and hence} \quad (\C[x]^o)^* = \prod_{\alpha \in \mathbb{A}^1_{\C}} \C[[x-\alpha]] \]
the natural inclusion $\C[x] \rInto (\C[x]^o)^*$ sending a polynomial to its Taylor-series expansion in every point $\alpha \in \mathbb{A}^1_{\C}$.

An intermediate step in arriving at $\mathbb{F}_1$-geometry would be to extend this complex coalgebra approach to integral schemes $\wis{Spec}(R)$ where $R$ is a finitely generated $\Z$-algebra, without additive torsion. In \cite{LeBruyn82} it was shown that in this case we still have Kostant duality, which asserts that for all $\Z$-algebras $R$ and all $\Z$-coalgebras $C$ there is a natural one-to-one correspondence
\[
\wis{alg}_{\Z}(R,C^*) \leftrightarrow \wis{coalg}_{\Z}(C,R^o) \]
if we take as the modified dual coalgebra $R^o$ the set of all $g \in R^* = Hom_{\Z}(R,\Z)$ such that $Ker(g)$ contains a twosided ideal $I \triangleleft R$ such that $R/I$ is a finitely generated projective $\Z$-module. 

The crucial difference with the complex case is that now that the relevant ideals $I$ no longer need to be comaximal and that there is no longer decomposition of the dual coalgebra. In our example when $R= \Z[x]$ the relevant ideals are those generated by a monic polynomial $f$ which can be decomposed in irreducible monic polynomials $f = g_1^{e_1} \hdots g_k^{e_k}$. But, as it may happen that $(g_i,g_j) \not= \Z[x]$ we have
\[
\frac{\Z[x]}{(f)} \not= \frac{\Z[x]}{(g_1)^{e_1}} \oplus \hdots \oplus \frac{\Z[x]}{(g_k)^{e_k}} \]
and we can no longer decompose the dual coalgebra $\Z[x]^o$ over the codimension one points $\mathbb{V}(g_i)$. Hence, we must recourse to describe the dual coalgebra as a direct limit
\[
\Z[x]^o = \underset{\rightarrow}{lim}~(\frac{\Z[x]}{(f)})^* \]
where the limit is considered with respect to divisibility of monic polynomials, as there are canonical inclusions of $\Z$-coalgebras
\[
(\frac{\Z[x]}{(f)})^* \rInto (\frac{\Z[x]}{(g)})^* \qquad \text{whenever} \qquad f | g \]
But then, also the $\Z$-algebra of distribution must be described as an inverse limit and we have a canonical ringmorphism
\[
\Z[x] \rInto (\Z[x]^o)^* = \underset{\leftarrow}{lim}~\frac{\Z[x]}{(f)} \]
Finally, to get at $\F_1$-geometry via this coalgebra approach we start with a $\lambda$-ring $R$ and define the $\lambda$-dual coalgebra
\[
R^o_{\lambda} = \{ g \in R^*~|~\exists I \subset \wis{ker}(g)~|~R/I~\text{$\lambda$-ring finite over $\Z$}~\} \]
which is indeed a coalgebra as the tensor product of $\lambda$-rings is again a $\lambda$-ring, or specialize even further to the {\em cyclotomic dual coalgebra} \index{cyclotomic dual coalgebra} $R^o_{cycl}$ on the sub-coalgebra of $R^o_{\lambda}$ spanned by the  maps $g$ having in their kernel an ideal $I$ such that $R/I \simeq \Z[x]/(\phi_1) \times \hdots \times \Z[x]/(\phi_k)$ where the $\phi_i$ are products of cyclotomic polynomials $\Phi_n(x)$. 

For example, the (cyclotomic) coalgebra representing $\PP^1_{\F_1}$ would then be
\[
C_{\PP^1_{\F_1}} = \Z[t_x] \oplus \underset{\rightarrow}{lim}~(\frac{\Z[x]}{(\phi)})^* \oplus \Z[t_{x^{-1}}] \]
where the $\phi$ run in the multiplicative system generated by the cyclotomic polynomials $\Phi_n(x)$ with $n \in \N_0$, the other two factors which are the enveloping coalgebras of the one-dimensional Lie algebra correspond to the points $[0]$ and $[\infty]$. Its corresponding algebra of distributions is then
\[
\Z[[x]] \oplus \widehat{\Z[x]}_{Hab} \oplus \Z[[x^{-1}]] \]
where $\widehat{\Z[x]}_{Hab}$ is the {\em Habiro ring} \index{Habiro ring} or the {\em cyclotomic completion} \index{cyclotomic completion} of $\Z[x]$ introduced and studied by Kazuo Habiro in \cite{Habiro02}. 

The Habiro ring is the straightforward generalization along the geometric axis of the profinite integers $\hat{Z}$ along the arithmetic axis. For we can write it as
\[
\widehat{\Z[x]}_{Hab} = \underset{\leftarrow}{lim}~\frac{\Z[x,x^{-1}]}{ [n!]_x} \qquad \text{with} \qquad [n!]_x = (x^n-1)(x^{n-1}-1) \hdots (x-1) \]
Its elements have a unique description as formal Laurent polynomials over $\Z$ of the form
\[
\sum_{n=0}^{\infty} a_n(x)[n!]_x \in \Z[[x,x^{-1}]] \qquad \text{with} \qquad deg(a_n(x)) < n \]
and hence can be evaluated at every root of unity (but possibly nowhere else). 

Some of its elements had been discovered before. For example Maxim Kontsevich observed in his investigations on Feynman integrals that the formal power series $\sum_{n=0}^{\infty} (-1)^n [n!]_x$ is defined in all roots of unity and Don Zagier subsequently proved the hilarious identity
\[
\sum_{n=0}^{\infty} (-1)^n [n!]_x = -\frac{1}{2} \sum_{n=1}^{\infty} n \chi(n) x^{\frac{(n^2-1)}{24}} \]
where $\chi$ is the quadratic character of conductor $12$, whereas the functions on both sides never makes sense simultaneously. The right hand side converges only within the unit circle, but still if one approaches a root of unity radially, the limit of the function values on the right coincide with the value on the left. Such functions are said to 'leak through' roots of unity. 

The Habiro topology was introduced to describe the properties of the Habiro ring $\widehat{\Z[x]}_{Hab}$. For example if $U$ is an infinite set of roots of unity having $\alpha \in \pmb{\mu}_{\infty}$ as a limit-point, meaning that $U$ contains infinitely many elements adjacent to $\alpha$, then if $f \in \widehat{\Z[x]}_{Hab}$ evaluates to zero in all roots $\beta \in U$ then $f = 0$. For any subset $S \subset \N_0$ define the completion
\[
\Z[x,x^{-1}]^S = \underset{\underset{\phi \in \Phi^*_S}{\leftarrow}}{lim}~\frac{\Z[x,x^{-1}]}{(\phi)} \]
where $\Phi^*_S$ is the multiplicative set of all monic polynomials generated by all cyclotomic polynomials $\Phi_n(x)$ for $n \in S$. Among the many precise results proved in \cite{Habiro02} we mention these two
\begin{itemize}
\item{If $S' \subset S$ with the property that every component of $S$ with respect to the nearness relation contains an element of $S'$, then the natural map between the completions in an inclusion
\[
\rho^S_{S'}~:~\Z[x,x^{-1}]^S \rInto \Z[x,x^{-1}]^{S'} \]}
\item{If $S$ is a {\em saturated subset} \index{saturated subset} of $\N_0$, which means that for all $n \in S$ also its divisor set $\langle n \rangle = \{ m | n \}$ is contained in $S$, then we have
\[
\Z[x,x^{-1}]^S = \bigcap_{n \in S} \Z[x,x^{-1}]^{\langle n \rangle} = \bigcap_{n \in S} \widehat{\Z[x,x^{-1}]}_{(x^n-1)} \]
where the right-hand terms are the $I$-adic completions of $\Z[x,x^{-1}]$ with respect to the ideals $I=(x^n-1)$.}
\end{itemize}

\section{Conway's big picture}

In \cite{Conway96}, John H. Conway investigates $\Q$-projectivity classes of lattices commensurable with the standard $2$-dimensional lattice $L_1=\langle \pmb{e}_1,\pmb{e}_2 \rangle = \Z \pmb{e}_1 + \Z \pmb{e}_2$ and he shows that any such lattice has a unique form
\[
L_{M,\frac{g}{h}} = \langle M \pmb{e}_1 + \frac{g}{h} \pmb{e}_2, \pmb{e}_2 \rangle \]
with rational numbers $M > 0$ and $0 \leq \frac{g}{h} < 1$. Lattices $L_{M,0}=L_M$ are called {\em number-like} \index{number-like lattice} and if, in addition, $M \in \N_0$ we just call them {\em numbers} \index{number lattice}.

We will now define a metric on the set of (equivalence classes) of lattices. For two lattices $L=L_{M,\frac{g}{h}}$ and $L'=L_{N,\frac{i}{j}}$ consider the matrix
\[
D_{L L'} = \begin{bmatrix} M & \frac{g}{h} \\ 0 & 1 \end{bmatrix} . \begin{bmatrix} N & \frac{i}{j} \\ 0 & 1 \end{bmatrix}^{-1} \]
and let $\alpha$ be the smallest positive rational number such that all entries of the matrix $\alpha D_{L L'}$ are integers. The {\em hyper-distance} \index{hyper-distance} between the lattices $L$ and $L'$ is then defined to be the integer
\[
\delta(L,L') = det(\alpha.D_{L L'}) \in \Z \]
One can show that the hyper-distance is symmetric and that $log(\delta(L.L'))$ is an ordinary metric on the projectivity classes of commensurable lattices. 

Conway's {\em big picture}\index{big picture} $\mathbb{B}$ is the graph with vertices the (classes of) lattices commensurable with $L_1$ and there is ad edge between the lattices $L$ and $L'$ if and only if $\delta(L,L')=p$, a prime number.  Conway shows that the sub-graph consisting of all lattices whose hyper-distance to $L_1$ is a power of $p$ is the infinite $p$-adic tree $T_p$, that is a $(p+1)$-valent tree, since for example the $p$-neighbors of $L_1$ are the lattices $L_p$ and $L_{\frac{1}{p},\frac{k}{p}}$ for $0 \leq k < p$. It must be a tree as the first step of the shortest path to $L_1$ from $L_{p^j}$ must be to $L_{p^{j-1}}$ as the other possibilities $L_{p^{j+1}}$ and $L_{p^j\frac{k}{p}}$ all have hyperdistance $p^{j+1}$ from $L_1$. Further he shows that the big picture is the product $\mathbb{B}=\ast_p T_p$. Here's part of the $2$-tree
\[
\begin{xymatrix}@=0.6em{
& & \frac{1}{8},\frac{5}{8} \ar@{-}[rd] & & \frac{1}{8},\frac{1}{8} \ar@{-}[ld] &  & \frac{1}{4},\frac{1}{8} \ar@{-}[rd] &  & \frac{1}{4},\frac{3}{8} \ar@{-}[ld] &  &  \\
& & & \frac{1}{4},\frac{1}{4} \ar@{-}[dd] & & & & \frac{1}{2},\frac{1}{4} \ar@{-}[dd] & & & \\
\frac{1}{8},\frac{3}{8} \ar@{-}[rd] & & & & & & & & & & \frac{1}{4},\frac{3}{8} \ar@{-}[ld] \\
& \frac{1}{4},\frac{3}{4} \ar@{-}[rr] & & \frac{1}{2},\frac{1}{2} & & & & 1\frac{1}{2} \ar@{-}[rr] & & \frac{1}{2},\frac{3}{4} & \\
\frac{1}{8},\frac{7}{8} \ar@{-}[ru] & & & & & & & & &  & \frac{1}{4},\frac{5}{8} \ar@{-}[lu] \\
& & & & 1 \ar@{-}[luu] \ar@[red]@{-}[rr] & & 2 \ar@{-}[ruu] & & & & \\
\frac{1}{8},\frac{1}{4} \ar@{-}[rd] & & & & & & & & & & 1\frac{3}{4} \ar@{-}[ld] \\
& \frac{1}{4},\frac{1}{2} \ar@{-}[rr] & & \frac{1}{2} \ar@[blue]@{-}[ruu] & & & & 4 \ar@[red]@{-}[luu] \ar@{-}[rr] & & 2\frac{1}{2} & \\
\frac{1}{8},\frac{3}{4} \ar@{-}[ru] & & & & & & & & & & 1\frac{1}{4} \ar@{-}[lu] \\
& & & \frac{1}{4} \ar@[blue]@{-}[uu] & & & & 8 \ar@[red]@{-}[uu] & & & \\
& & \frac{1}{8},\frac{1}{2} \ar@{-}[ru] & & \frac{1}{8} \ar@[blue]@{-}[lu] & & 16 \ar@[red]@{-}[ru]  & & 4\frac{1}{2} \ar@{-}[lu] & &}
\end{xymatrix}
\]
Sometimes it is helpful to choose another normalization for the lattice $L$ by swapping the vectors $\pmb{e}_1$ and $\pmb{e}_2$. Let $\pmb{v}_1 = M \pmb{e}_1 + \tfrac{g}{h} \pmb{e}_2$ and $\pmb{v}_2 = \pmb{e}_2$ be the standard generators of $L=L_{M,\frac{g}{h}}$ then $L$ is also generated by the vectors
\[
h \pmb{v}_1 - g \pmb{v}_2 = hM \pmb{e}_1 \quad \text{and} \quad g' \pmb{v}_1 - h' \pmb{v}_2 = g' M \pmb{e_1} + \frac{1}{h} \pmb{e}_2 \]
where $g',h' \in \Z$ such that $gg'-hh'=1$. Dividing by $hM$ we get the {\em reversed normalized} form for $L_{M,\frac{g}{h}}$
\[
L = \langle \frac{1}{h^2 M} \pmb{e}_2 + \frac{g'}{h} \pmb{e}_1,\pmb{e}_1 \rangle \]
So we get an involution on the vertices of the big picture
\[
(M,\frac{g}{h}) \leftrightarrow (\frac{1}{h^2M},\frac{g'}{h}) \]
where $g'$ is the inverse of $g$ modulo $h$. 

The vertices of the big picture correspond to couples $(M,\frac{g}{h})$ so are elements of $\Q_0 \times \Q/\Z$ and we can identify each of the factors $\Q/\Z$ (additively) with $\pmb{\mu}_{\infty}$ (multiplicatively). One quickly verifies that the hyperdistance
\[
\delta(L_M,L_{M,\frac{g}{h}}) = h^2 \]
So the cyclic subgroup $\pmb{\mu}_n$ corresponding to $M$ is contained in a ball $B(L_M,n^2)$ around the lattice $L_M$ with hyper-distance $n^2$. In particular,   the non-identity elements of cyclic group $\pmb{\mu}_p$ for $p$ a prime number have hyper-distance $p^2$ from $L_M$ and are the $p-1$ vertices in the $p$-tree connected to $L_{pM}$. 

The lattices $L_n$ with $n \in \N_0$ form the big cell in this picture, which is the product of graphs of type $A_{\infty}^+$, one for each prime number $p$
\[
A_{\infty}~:~\xymatrix{1 \ar[r] & p \ar[r] & p^2 \ar[r] & \hdots & \hdots \ar[r] & p^k \ar[r]  & \hdots} \] and can be identified with $\N_0^{\times} = \PP^1_{\F_1} - \{ [0],[\infty] \}$. But then, we can extend the Habiro topology to Conway's big picture by calling two lattices related
\[
L \sim L' \qquad \Leftrightarrow \qquad \delta(L,L')=p^a \]
if their hyper-distance is a pure prime-power. An open set is then a subset $U$ of vertices of the big picture having the property
that for each $L \in U$ the set $\{ L' \sim L~|~L' \notin U \}$ is finite. Clearly, the restriction of the Habiro topology on the big picture restricts to the usual Habiro topology on the big cell $\N_0^{\times}$.

$\Z \mathbb{B}$, the free $\Z$-module on the vertices of $\mathbb{B}$, is the playground of several operations on $\mathbb{B}$. Some well-known classical ones such as the {\em Hecke operators} \index{Hecke operators} $T_n$ which takes the vertex representing the lattice $L$ to the sum of all vertices corresponding to lattices $L'$ with $\delta(L,L')=n$, that is $T_p$ replaces the center of each ball of hyperradius $n$ by its periphery. These Hecke operators clearly satisfy for $a > 1$ the relation
\[
T_{p} \circ T_{p^a} = p T_{p^{a-1}} + T_{p^{a+1}} \]
as the left-hand takes a vertex to the sum of all neighbors of vertices at hyperdistance $p^a$ from it, but in this sum each vertex of hyperdistance $p^{a-1}$ occurs $p$ times and each point of hyperdistance $p^{a+1}$ just once, giving the right-hand side. 

More operators come from the action of a certain noncommutative algebra on $\Z \mathbb{B}$, the {\em Bost-Connes algebra} \index{Bost-Connes algebra} $\Lambda$, see for example \cite{FUN}. If $\Lambda$ hadn't been constructed years before (in \cite{BC95}) it would have arisen naturally from $\F_1$-geometry by a construction well-known in noncommutative algebraic geometry. 

If $X$ is an affine $\C$-variety with a linear action by a finite group $G$, then the coordinate ring of its quotient variety $\C[X/G]=\C[X]^G$ is Morita equivalent (that is, have equivalent module categories) to the skew group algebra $\C[X] \ast G$ which as a $\C$-vectorspace is the group-algebra $\C[X] G$ but with multiplication induced by $f.e_g = e_g.\phi_g(f)$. That is, one way to get at the descended algebra $\C[X]^G$ is by considering the noncommutative skew group algebra $\C[X] \ast G$. 

In Borger's proposal for $\F_1$-geometry this approach may be very helpful as an $\F_1$-algebra is a $\Z$-algebra $R$ together with descent data given by the action of the monoid $\N_0^{\times}$ by the endomorphisms $\{ \Psi^n \}$. But now we cannot directly construct the invariant algebra $R^{\N_0^{\times}}$ (which would be our elusive $\F_1$-algebra) but we can still construct the skew-monoid algebra $R \ast \N_0^{\times}$ which, as before,  as $\Z$-module is $R G = \oplus_{n \in \N_0^{\times}} R e_n$ but with noncommutative multiplication induced by $r.e_n = e_n.\Psi_n(r)$.

For example, let us try to understand the algebraic closure $\overline{\F_1}$ by considering the associated skew monoid algebra. The $\lambda$-algebra corresponding to $\overline{\F_1}$ is the group-algebra $\Z[\pmb{\mu}_{\infty}]$ with Frobenius-lifts $\Psi^n$ induced by sending a root of unity $\omega$ to $\omega^n$. If we write the group-law additively instead of multiplicative we get the group-algebra $\Z[\Q/\Z]$ with $\Psi^n(e(\frac{g}{h})) = e(\frac{ng}{h}~\wis{mod}~1)$. The corresponding skew-monoid algebra is then
\[
\Z[\Q/\Z] \ast \N_0^{\times} = \bigoplus_{n \in \N_0^{\times}} \Z[\Q/\Z] e_n \quad \text{with} \quad  e(\frac{g}{h}).e_n = e_n.\Psi^n(e(\frac{g}{h})) \]
Noncommutative algebraic geometers would then study properties of this ring to get insight into $\overline{\F_1}$. Noncommutative {\em differential} geometers however work with $\ast$-algebras, so they have to construct the minimal $\ast$-algebra generated by $\Z[\Q/\Z] \ast \N_0^{\times}$ and therefore consider the algebra
\[
\Z[\Q/\Z] \bowtie \N_0^{\times} = \bigoplus_{m,n \in \N_0^{\times}, (m,n)=1} e^*_m \Z[\Q/\Z] e_n \]
in which the generators $e^*_m,e_n$ and $e(\frac{g}{h})$ satisfy the following multiplication rules
\[
\begin{cases}
e_n.e(\frac{g}{h}).e^*_n = \rho_n(e(\frac{g}{h})) \\
e^*_n.e(\frac{g}{h}) = \Psi^n(e(\frac{g}{h})) e^*_n \\
e(\frac{g}{h}).e_n = e_n.\Psi^n(e(\frac{g}{h}) \\
e_n.e_m = e_{nm} \\
e^*_n.e^*_m = e^*_{nm} \\
e^*_n.e_n = n \\
e_n.e^*_m = e^*_m.e_n~\text{if $(m,n)=1$}
\end{cases}
\]
where $\rho_n(e(\frac{g}{h})) = \sum_{n.\frac{i}{j} = \frac{g}{h}} e(\frac{i}{j})$. This algebra $\Z[\Q/\Z] \bowtie \N_0^{\times}$ is then the (integral version of the) Bost-Connes algebra $\Lambda$ constructed in \cite{FUN} where it is also shown that there is an action of $\Lambda$ and $\Z \mathbb{B}$ given by the rules
\begin{itemize}
\item{$e_n.L_{\frac{c}{d},\frac{g}{h}} = L_{\frac{nc}{d},\rho^m(\frac{g}{h})}$ where $m=(n,d)$}
\item{$e^*_n.L_{\frac{c}{d},\frac{g}{h}} = (n,c) L_{\frac{c}{nd},\Psi^{\frac{n}{m}}(\frac{g}{h})}$ where $m=(n,c)$}
\item{$e(\frac{a}{b}).L_{\frac{c}{d},\frac{g}{h}} = L_{\frac{c}{d},\Psi^c(\frac{a}{b}) \frac{g}{h}}$}
\end{itemize}
with $\rho_n$ and $\Psi^n$ defined on $\frac{g}{h}$ as they were defined before on $e(\frac{g}{h})$. So far, we have identified $\PP^1_{\F_1}$ (equipped with the adjacency relation among its schematic points) with the big cell in the Conway picture. It is believed that this bigger picture will play an ever increasing role of importance in future developments in $\F_1$-geometry and will illuminate surprise appearances of the Bost-Connes algebra $\Lambda$ as  a generalized symmetry on geometric $\F_1$-objects, see for example \cite{ConnesConsani11}.

\section{Exotic topology on $\wis{Spec}(\Z)$}

Now that we have a formal definition of $\PP^1_{\F_1}$ let's try to make sense of the ultimate question in $\F_1$-geometry : what (if any) geometric object is $\wis{Spec}(\Z)$ over $\F_1$? Again, we will start with an intuitive proposal due to A. L. Smirnov \cite{Smirnov92} and later try to formalize it using $\lambda$-rings.

Smirnov proposes to take as the set of {\em schematic points} of $\overline{\wis{Spec}(\mathbb{Z})}$ the set
\[
\{ [2],[3],[5],[7],[11],[13],[17],\hdots \} \cup \{ [ \infty ] \} \]
of all prime numbers together with a point at infinity. The {\em degrees} of these schematic points are then taken to be
\[
\wis{deg}([p]) = log(p) \qquad \text{and} \qquad \wis{deg}([\infty ]) = 1 \]
To understand this proposal, recall that if $C$ is a smooth projective curve over $\mathbb{F}_p$, then a schematic point $P \in C$ corresponds to a discrete valuation ring $\Oscr_P$ in the function field $\mathbb{F}_p(C)$ with maximal ideal $\frak{m}_P = (t_P)$ where $t_P$ is a uniformizing parameter. The degree of the schematic point $P \in C$, $\wis{deg}(P)$, is defined to be $n$ if and only if $\Oscr_P/\frak{m}_P = \mathbb{F}_{p^n}$.
A rational function $f \in \mathbb{F}_1(C)$ is said to be {\em regular} in $P$ if and only if $f \in \Oscr_P$ and the {\em order} of $f$ in $P$ is the valuation of $f$, that is, $\wis{ord}_P(f) = k$ if and only if $f \in \frak{m}_P^k - \frak{m}_P^{k+1}$ for a unique $k \in \mathbb{Z}$. The {\em divisor} of the rational function $f \in \mathbb{F}_p(C)$
\[
\wis{div}(f) = \sum_{P \in C} \wis{ord}_P(f) P \]
then has degree zero, that is
\[
0 = \wis{deg}(\wis{div}(f)) = \sum_{P \in C} \wis{ord}_P(f) \wis{deg}(P) \]
By analogy, we may take the schematic points of $\overline{\wis{Spec}(\Z)}$ to be the different discrete valuation rings in the corresponding 'function field' $\mathbb{Q}$. By Ostrovski's theorem they are either the {\em $p$-adic valuations}
\[
v_p(q) = n \qquad \text{if} \qquad q =p^n \frac{r}{s}~\text{and}~p \nmid r.s \]
for every prime number $p$ together with the {\em real valuation}
\[
v_{\infty}(q) = - log( |q| ) \]
the minus sign arrises because of the convention that $v_{\infty}(0) = \infty$. But then, if $q = \pm \frac{p_1^{e_1} \hdots p_r^{e_r}}{q_1^{f_1} \hdots q_s^{e_s}}$ its corresponding divisor must be
\[
\wis{div}(q) = \sum_{i=1}^r e_i [p_i] - \sum_{j=1}^s f_j [q_j] - log( | q |) [\infty] \]
The proposal for the degrees of the schematic points of $\overline{\wis{Spec}(\Z)}$ is then the only one possible (up to a common multiple) such that the degrees of all these principal divisors are equal to zero.
Any non-constant rational function $f \in \mathbb{F}_p(C)$ determines a {\em cover} map $f~:~C \rOnto \mathbb{P}^1_{\mathbb{F}_p}$. Smirnov defines as the constant rational numbers the intersection $\Q \cap \overline{\F_1} = \{ 0 \} \cup \{ 1,-1 \} = \F_{1^2}$.  Therefore, we would expect by analogy a rational numbers $q= \tfrac{a}{b} \in \mathbb{Q}$ with $(a,b)=1$ to determine a cover
\[
q~:~\overline{\wis{Spec}(\Z)} \rOnto \mathbb{P}^1_{\mathbb{F}_1} \]
Smirnov's proposal in \cite{Smirnov92} is to {\em define} a map by
\[
[p] \mapsto \begin{cases} [0]~\quad \text{if $p | a$} \\
[\infty]~\quad \text{if $p | b$} \\
[n]~\quad \text{if $p \nmid ab$ and $\overline{a} \overline{b}^{-1}$ has order $n$ in $\F_p^*$}
\end{cases} \]
and by sending
\[
[\infty ] \mapsto \begin{cases} [0]~\quad \text{if $a < b$} \\ [\infty]~\quad \text{if $a > b$} \end{cases} \]

To motivate this definition let us again look at the function-field case. Any rational function $f \in \F_p(C)$ defines a map between the geometric points
\[
C(\overline{\F}_p) \rTo \mathbb{P}^1_{\F_p}(\overline{\F}_p) \qquad P \mapsto \begin{cases} [f(P):1] & \text{if $f \in \Oscr_P$} \\ [\infty] &  \text{if $f \notin \Oscr_P$} \end{cases} \]
with $f(P) = \overline{f} \in \Oscr_P/\frak{m}_P = \overline{\F}_p$. Because $f \in \F_p(C)$ we have for all $P \in C(\overline{\F}_p)$ and all $\sigma \in Gal(\overline{\F}_p/\F_p)$ that $\sigma(f(P)) = f(\sigma(P))$ and hence this map induces a map between the schematic points $C \rTo \mathbb{P}^1_{\F_p}$ sending a schematic point (a Galois orbit of a geometric point $P$) to the Frobenius-orbit of the root of unity $f(P)$ (or its corresponding monic irreducible polynomial in $\F_p[x]$).
Returning to the above map given by a rational number $q = \frac{a}{b}$ it is clear that $q([p])=[0]$ for all prime factors $p$ of $a$ and that $q([p])=[\infty]$ for all factors of $b$. To understand the other images note that if $\overline{a} \overline{b}^{-1}$ has order $n$ in $\F_p^*$, there exists a prime ideal $P$ in the ring of cyclotomic integers $\Z[\epsilon]$ (for $\epsilon$ a primitive $n$-th root of unity) lying over $(p)$, that is $P \cap \Z = (p)$ with corresponding discrete valuation ring $\Oscr_P$ such that
\[
\frac{a}{b} - \epsilon \in P \Oscr_P = \frak{m}_P \]
and therefore $\frac{a}{b}(P) = \epsilon(P)$, explaining why the schematic point $[p]$ is send to the Galois orbit of $\epsilon$ which is precisely the schematic point $[n]$ of $\mathbb{P}^1_{\F_1}$.

In the function-field case we have for every non-constant rational function $f \in \F_p(C) - \F_p$ that the map $C \rOnto \mathbb{P}^1_{F_p}$ is surjective with finite fibers. Let us first verify finiteness for the map $q = \frac{a}{b}$, that is, for every $[n]$ we must show that there are only finitely many primes $p$ for which
\[
(\frac{\overline{a}}{\overline{b}})^n = 1 \quad \text{in}~\F_p^* \]
This is clearly equivalent to $p | a^n-b^n$ {\em and} $p \nmid a^m-b^m$ for all $m < n$, so $q^{-1}([n])$ is a subset of the finite number of prime factors of $a^n-b^n$.
Surjectivity of the map $q$ is less clear as there seems to be no reason why there should be always a prime factor of $a^n-b^n$ not dividing the number $a^m-b^m$ for all $m < n$. In fact, surjectivity is not always true. 
\[
\includegraphics[width=12cm]{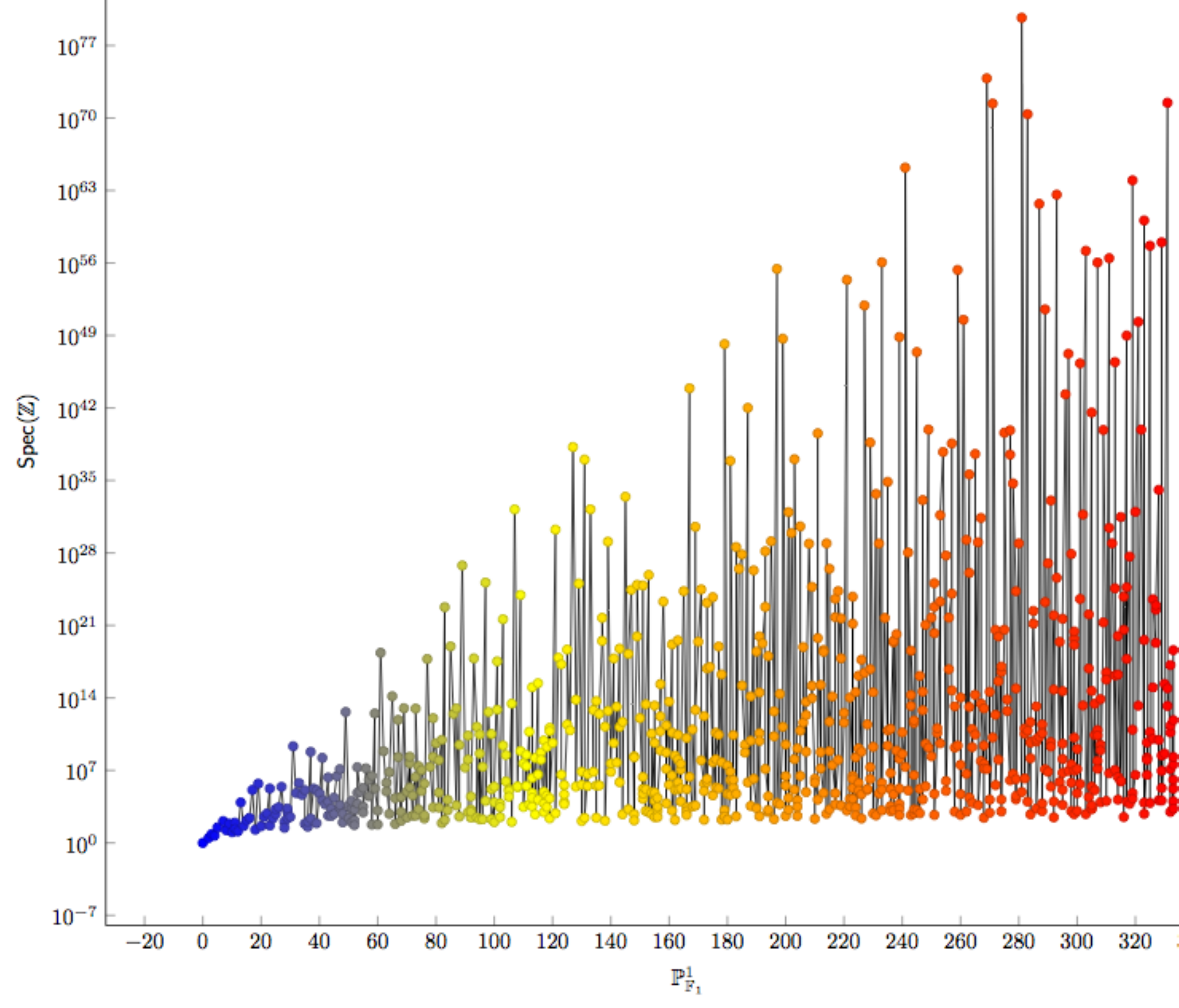} \]
For example, the map $q=\frac{2}{1}$ has no prime mapping to $[6]$.
The above picture gives a portion of the graph of the map $2$ in the Smirnov-plane $\mathbb{P}^1_{\F_1} \times \wis{Spec}(\Z)$, where we have used a logarithmic scale on the prime-number axis and determined the full fibers of all $[n] \in \PP^1_{\F_1}$ for $n < 330$.
 The points on the 'diagonal' are the first few Mersenne prime numbers, that is primes $p$ such that $M_p=2^p-1$ is again a prime number.

Perhaps surprisingly we can determine all rational numbers $q$ for which the map $\overline{\wis{Spec}(\Z)} \rTo \mathbb{P}^1_{\F_1}$ fails to be surjective, as well the schematic points $[n]$ of $\mathbb{P}^1_{\F_1}$ for which $q^{-1}([n]) = \emptyset$. The crucial result needed is {\em Zsigmondy's theorem}\index{Zsigmondy's theorem} \cite{Zsig}.
Consider positive integers $1 \leq b < a$ with $(a,b)=1$. Then, for every $n > 1$ there exist prime numbers $p | a^n-b^n$ such that $p \nmid a^m-b^m$ for all $m < n$ unless we are in one of the following two cases:
\begin{enumerate}
\item{$a=2, b=1$ and $n=6$, or}
\item{$a+b=2^k$ and $n=2$}
\end{enumerate}

Smirnov's interest in these maps is that the ABC-conjecture would follow provided one can prove a suitable analogue of the{\em  Riemann-Hurwitz formula} \index{Riemann-Hurwitz formula} for the maps $q$. Recall that if $f~:~C \rOnto \PP^1_k$ is a non-constant cover from a smooth projective curve $C$ over a field $k$, then the Riemann-Hurwitz formula asserts that
\[
2 g_C - 2 \geq -2 deg(f) + \sum_{scheme} (e_f(P)-1) deg(P) \]
where $g_C$ is the genus of $C$ and $e_f(P)$ is the ramification of a schematic point $P \in C$ of degree $deg)P)$. If we define the {\em defect} \index{defect} $\delta_P$ of a schematic point $P \in C$ to be the number $\delta_P = \tfrac{(e_f(P)-1)deg(P)}{deg)f)} \geq 0$, then the Riemann-Hurwitz formula can be rewritten as
\[
\sum_{scheme} \delta_P \leq 2 - \frac{2 - 2g_C}{deg(f)} \]
and we note that this inequality still holds if we sum over a sub-selection of the schematic points $P \in C$. Again, we want to define via analogy the {\em ramification index } \index{ramification index} $e_q(p)$ and the {\em arithmetic defect} \index{arithmetic defect} $\delta(p)$ for any prime number $p$ with respect to a cover $q~:~\overline{\wis{Spec}(\Z)} \rOnto \PP^1_{\F_1}$. If $q = \tfrac{a}{b}$, then Smirnov proposes to take for $e_q(p)$ the largest power of $p$ dividing $a$ (provided $p \in q^{-1}([0])$), the largest power dividing $b$ (provided $p \in q^{-1}([\infty])$) and if $p \in q^{-1}([n])$ to take as $e_q(p) = k$ if $p^k$ is the largest power dividing $a^n-b^n$. With this definition of the ramification index, het then proposes to define the arithmetic defect
\[
\delta(p) = \frac{(e_q(p)-1) log(p)}{log(a)} \]
which coincides with the classical definition (given our proposal for the degree of $p$) provided we define the degree of the map $q$ to be $log(a)$. Let us try to motivate this proposal in the case when $a >> b$. Take $[p] \in \PP^1_{\F_1}$ for $p$ a prime number, then the divisor of $q^{-1}([p]) = \sum_i n_i q_i$ if $\prod_i q_i^{n_i} = a^p-b^p$, whence $\sum_i n_i log(q_i) \approx p.log(a)$, and as we took $deg([p])=\phi(p)=p-1$ it follows that indeed 
\[
\frac{\sum_i n_i deg(q_i)}{deg([p])} \approx log(a) \]
For any schematic point $[d] \in \PP^1_{\F_1}$ let us define the defect of $[d]$ to be
\[
\delta([d]) = \sum_{p \in q^{-1}([d])} \delta(p) \]
Now, if $a=p_1^{e_1}\hdots p_k^{e_k}$ and $b = q_1^{f_1} \hdots q_s^{f_s}$ we define $a_0=p_1 \hdots p_k$ and $b_0 = q_1 \hdots q_s$ and $a_1=a/a_0$, $b_1=b/b_0$. Then, with the above proposals it is easy to work out that
\[
\delta([0]) = \frac{log(a_1)}{log(a)}, \qquad \delta([\infty]) = \frac{log(b_1) + log(q)-1}{log(a)}, \qquad \delta([1]) = \frac{log((a-b)_1)}{log(a)} \]
as $q^{-1}([1]) = \{ p~|~a-b \}$ and where in the middle term $log(q)-1$ is the contribution of $\infty$ in the defect. If we would be able to prove a variant of the Riemann-Hurwitz formula in $\F_1$-geometry for all covering maps $q~:~\overline{\wis{Spec}(\Z)} \rOnto \PP^1_{\F_1}$ and if we assume the constant $\gamma = 2 g_{\overline{\wis{Spec}(\Z)}} -2 \geq 0$, then it would follow that (limiting to points lying in the fibers of $[0],[1]$ and $[\infty]$ that
$\delta([0]) + \delta([1]) + \delta([\infty]) =$
\[
\frac{1}{log(a)}(log(a_1)+log((a-b)_1)+log(b_1)+log(a)-log(b)-1) \leq 2 + \frac{\gamma}{log(a)} \]
Now, let's turn to the ABC-conjecture \index{ABC-conjecture}. Suppose $A+B=C$ with $(A,B,C)=1$ and take $a=C$ and $b=min(A,B)$ and consider the cover $q = \tfrac{a}{b}~:~\overline{\wis{Spec}(\Z)} \rOnto \PP^1_{\F_1}$ and clearly we have $a-b \geq \tfrac{a}{2}$. Then, we can rewrite the above inequality as
\[
1 \leq \frac{log(a_0.b_0.(a-b)_0}{log(a)} + \frac{log(C')}{log(a)} \]
where $log(C')=\gamma+log(2)+1$, but then $a \leq C'(a_0.b_0.(a-b)_0)$ or in other words that
\[
C \leq C'(rad(A.B.C)) \]
which is a (too strong) formulation of the ABC-conjecture.

Now that we have a family of non-constant cover-maps for all $q \in \Q$
\[
q~:~\wis{Spec}(\Z) \rTo \PP^1_{\F_1} \]
and have the Habiro topology on $\PP^1_{\F_1}$ we can define the {\em exotic topology} \index{exotic topology} on $\wis{Spec}(\Z)$ to be the coarsest topology with the property that all maps $q$ are continuous with respect to the Habiro topology on $\PP^1_{\F_1}$. As the covers are finite and the Habiro topology is finer than the cofinite topology we have that the exotic topology refines the usual, that is cofinite, topology on $\wis{Spec}(\Z)$. Again, this topology is no longer compact.

\section{Witt and Burnside rings}

Surprisingly, the forgetful functor from $f~:~\wis{rings}_{\lambda} \rTo \wis{rings}$ has a {\em right adjoint} (a left adjoint is the common situation), that is there is a functor
\[
w~:~\wis{rings} \rTo \wis{rings}_{\lambda} \quad \text{such that} \quad \wis{alg}_{\lambda}(A,w(B)) = \wis{alg}(f(A),B) \]
for all $\lambda$-rings $A$ and all rings $B$. We will recall the construction of this {\em 'witty functor'} \index{witty functor} (it is closely related to the functor of Big Witt vectors).

For any ring $A$ let $w(A) = 1+tA[[t]]$ be the set of all formal power series with coefficients in $A$ and with constant term equal to $1$. We will turn this set into a rings with addition $\oplus$ and multiplication $\otimes$ to distinguish these operations from the usual ones on the formal power series ring $A[[t]]$. The addition $\oplus$ on $w(A)$ will be the usual multiplication of formal power series, that is
\[
u(t) \oplus v(t) = u(t).v(t) \quad \text{and hence} \quad \boxed{0}=1 \quad \text{and} \quad \ominus u(t) = u(t)^{-1} \]
Multiplication is enforced by functoriality and the rule that for all $a,b \in A$ we demand that
\[
\frac{1}{1-at} \otimes \frac{1}{1-bt} = \frac{1}{1-abt} \quad \text{and hence} \quad \boxed{1} = \frac{1}{1-t}= 1+t+t^2+\hdots \]
What we mean by this is that for any $u(t) \in w(A)$ there exists unique $a_i \in A$ such that
\[
u(t) = \prod_{i=1}^{\infty} \frac{1}{1-a_it^i} \]
For each $n$, denote $\alpha_n = \sqrt[n]{a_n}$ and let $\zeta_n$ be a primitive $n$-th root of unity, so that for all $n$ we have that $1-a_nt^n = \prod_{i=0}^n (1-\zeta_n^i \alpha_n)$. But then, over the ring $A[\pmb{\mu}_{\infty}][\alpha_1,\alpha_2,\hdots]$ we can write $u(t)$ as
\[
u(t) = A_1 \oplus A_2 \oplus A_3 \oplus \hdots \quad \text{with} \quad A_n = \frac{1}{1-\alpha_nt} \oplus \frac{1}{1-\zeta_n\alpha_n t} \oplus \hdots \oplus \frac{1}{1-\zeta_n^{n-1}\alpha_n t} \]
If we similarly write the power series $v(t) = B_1 \oplus B_2 \oplus B_3 \oplus \hdots$, then the product must be
\[
u(t) \otimes v(t) = C_1 \oplus C_2 \oplus \hdots \quad \text{with} \quad C_{i+1} = \bigoplus_{j+k=i+1} A_j \otimes B_k \]
and by construction and symmetric function theory one verifies that the formal power series $u(t) \otimes v(t)$ has all its coefficients in $A$. In this way we see that $w(A)$ is a commutative ring with zero-element the constant power series $1$ and multiplicative unit the power series $1+t+t^2+\hdots$.
In addition, $w(A)$ becomes a $\lambda$-ring with the Frobenius lifts induced by the rule that
\[
\psi^p(\frac{1}{1-at}) = \frac{1}{1-a^nt} \]
and extended additively so that if $u(t) = A_1 \oplus A_2 \oplus \hdots $ then $\psi^p(u(t)) = \psi^p(A_1) \oplus \psi^p(A_2) \oplus \hdots $. The Frobenius lifts are also multiplicative by functoriality and the calculation that
\[
\psi^p(\frac{1}{1-at} \otimes \frac{1}{1-bt}) = \frac{1}{1-a^pb^pt} = \psi^p(\frac{1}{1-at}) \otimes \psi^p(\frac{1}{1-bt}) \]
Clearly, the endomorphisms $\psi^n, n \in \N$ commute with each other and $\psi^p$ is a Frobenius lift as
\[
(\frac{1}{1-a_1t} \oplus \hdots \oplus \frac{1}{1-a_kt})^{\otimes p} - (\frac{1}{1-a_1^pt} \oplus \hdots \oplus \frac{1}{1-a_k^pt}) \]
is divisible by $p$ by the binomial formula. There is an additional family of additive-group endomorphisms $V_n$ on $w(A)$, the {\em Verschiebungs-operators} \index{Verschiebung} which are defined by $V_n(s(t)) = s(t^n)$ and finally there is the $[n]$ operator which maps $s(t)$ to $s(t)^{n}=s(t) \oplus \hdots \oplus s(t)$ ($n$ times). These maps satisfy the following relations
\[
V_n \circ V_m = V_m \circ V_n \qquad \Psi^m \circ \Psi^n=\Psi^n \circ \Psi^m \qquad \Psi^n \circ V_n = [n] \]
and $\Psi^n \circ V_m = V_m \circ \Psi^n$ if $(m,n)=1$.
This witty construction is functorial because for any ring-morphism $\phi~:~A \rTo B$ we have a ringmorphism $\Phi~:~w(A) \rTo w(B)$ compatible with the Frobenius lifts, induced by the rule that
\[
\Phi(\frac{1}{1-at}) = \frac{1}{1-\phi(a)t} \]
which gives us that $\Phi(1+a_1t+a_2t^2+\hdots) = 1 + \phi(a_1)t+\phi(a_2)t^2+\hdots$. 
We now define maps $\gamma_n~:~w(A) \rTo A$ via the formula
\[
\frac{t u'}{u} = \sum_{n=1}^{\infty} \gamma_n(u) t^n \]
where we have used the {\em logarithmic derivative} \index{logarithmic derivative} $\tfrac{u'}{u}$ which transforms multiplication in addition. If we work this out for $u= \tfrac{1}{1-at}$ then $u' = \tfrac{a}{(1-at(^2}$ and hence
\[
\frac{t u'}{u} = at + a^2t^2+ a^3t^3 + \hdots \]
whence $\gamma_n(\tfrac{1}{1-at}) = a^n$ and is therefore multiplicative in $a$. Again, using functoriality it is then easy to conclude that all the maps $\gamma_n~:~w(A) \rTo A$ are in fact ringmorphisms.

If $A$ is in addition a $\lambda$-ring with commuting family of endomorphisms $\Psi^n$ generated by the Frobenius lifts, then there is a $\lambda$-ring morphism $\sigma_t$ making the diagram below commute
\[
\xymatrix{ & & & & & w(A) = 1+tA[[t]] \ar[d]^{(\gamma_1,\gamma_2,\gamma_3,\hdots)} \\
A \ar[rrrrru]^{\sigma_t} \ar[rrrrr]_{\Psi=(\Psi^1,\Psi^2,\Psi^3,\hdots)} & & & & & A^{\omega} = (A,A,A,\hdots)}
\]
where $\sigma_t$ is defined via the formula
\[
\sigma_t(a) = exp(\int \frac{1}{t} \sum_{n=1}^{\infty} \Psi^n(a) t^n) \]
Again, it is easy to verify that $\sigma_t(a+b) = \sigma_t(a) \oplus \sigma_t(b)$ and slightly more difficult to prove that $\sigma_t(a.b) = \sigma_t(a) \otimes \sigma_t(b)$, whence $\sigma_t$ is a ringmorphism and is compatible with the $\Psi^n$-endomorphisms, so is a $\lambda$-ringmorphism.

From these facts the right-adjointness of the witty functor with respect to the forgetful functor follows. If $A$ is a $\lambda$-ring and $\phi$ is a ringmorphism $f(A)=A \rTo B$, then we get a $\lambda$-ring morphism
\[
A \rTo^{\sigma_t} w(A) \rTo^{\Phi} w(B) \]
Conversely, a $\lambda$-ring morphism $A \rTo w(B)$ composed with the ringmorphism $\gamma_1~:~w(B) \rTo B$ gives a ringmorphism $A \rTo B$ and one verifies both constructions are each other inverses.

If one accepts Borger's proposal that $\mathbb{F}_1$-algebras are just $\lambda$-rings without additive torsion, where we interpret the commuting family of endomorphisms $\{ \Psi^n~:~n \in \N_0 \}$ as descent data from $\Z$ to $\mathbb{F}_1$, then the forgetful functor
\[
f = - \otimes_{\F_1} \Z~:~\wis{alg}_{\F_1} = \wis{rings}_{\lambda} \rTo \wis{rings} \]
that is, stripping off the descent data, can be interpreted as the base-extension functor from $\F_1$ to $\Z$. But then, as a right adjoint to base-extension, the witty functor $w$ can be interpreted as the Weil descent from $\Z$-rings to $\mathbb{F}_1$-algebras. Hence, we finally know what $\wis{Spec}(\Z)$ should be over the elusive field $\F_1$ : it must be the geometric object associated to the $\lambda$-algebra $w(\Z)$!

We will now make the connection between the construction of $w(A)$ and the more classical notion of the ring of Big Witt vectors $W(A)$. For much more detail we refer to the lecture notes of Michiel Hazewinkel \cite{Hazewinkel08}. Let us take $W(A) = A^{\omega} = (A,A,A,\hdots)$ and consider the diagram
\[
\xymatrix{ W(A) \ar[rrr]^{\gamma} & & & w(A) \ar[d]^{\tfrac{u'}{u}} \\
A^{\omega} \ar[u]^{=} \ar[rrr]_{=} & & & tA[[t]]}
\]
where $\gamma$ is the map sending $(a_1,a_2,a_3,\hdots)$ to $\prod_i \tfrac{1}{1-a_it^i}$ and which can be used to define a ringstructure on the {\em Big Witt vectors} \index{Big Witt vectors} $W(A)$ by transport of structure.

Before we describe the geometry, let us give a combinatorial interpretation of $w(\Z)$ due to Andreas Dress and Christian Siebeneicher \cite{DressSieb89}. 

Let $C=C_{\infty} = \langle c \rangle$ be the infinite cyclic group, written multiplicatively. A $C$-set $X$ is called {\em almost finite} if $X$ has no infinite orbits and if the number of orbits of size $n$ is finite for every $n \in \N_0$. A motivating example is the set of geometric points $X(\overline{\F_p})$ of an $\F_p$-variety $X$ on which $c$ acts as the Frobenius morphism. 

If $X$ and $Y$ are almost finite $C$-sets, then so are their disjoint union $X \sqcup Y$ and their Cartesian product $X \times Y$. These operations define an addition $+=\sqcup$ and multiplication $. = \times$ on the isomorphism classes $\hat{B}(C)$ of all almost-finite $C$-sets, called the {\em Burnside ring} \index{Burnside ring}. For any almost finite $C$-set $X$ define
\[
\phi_{C^n}(X) = \# \{ x \in X~:~c^n.x = x \} \]
that is, the number of elements lying in a $C$-orbit of sizes a divisor of $n$ and hence this number is finite. Moreover, the $\phi_{C^n}$ take disjoint unions (resp. products) to sums (resp. products) of the corresponding numbers and so all maps $\phi_{C^n}~:~\hat{B}(C) \rTo \Z$ are ringmorphisms. This gives us a collective ringmorphism
\[
\hat{\phi} = \prod_n \phi_{C^n}~:~\hat{B}(C) \rTo \Z^{\omega} = \wis{gh}(C) \]
where $\wis{gh}(C)$ is the {\em ghost ring} \index{ghost ring}, that is all maps $\N \rTo \Z$ with componentwise addition and multiplication. One verifies that $\hat{\phi}$ is injective but not surjective. We can extend the diagram of the previous section
\[
\xymatrix{
& Nr(\Z) \ar[d]^{itp} & \\
W(\Z) \ar[r]^{\tau} \ar[d]_{\Phi} & \hat{B}(C) \ar[r]^{s_t} \ar[d]_{\hat{\phi}} & w(\Z) \ar[d]^{L_z} \\
\prod_n \Z \ar[r]^{obv} & \wis{gh}(C) \ar[r]^{idn} & t \Z[[t]]}
\]
where $Nr(\Z)$ is called the {\em necklace algebra} \index{necklace algebra}, that is, the set $\Z^{\omega}$ with componentwise addition but multiplication defined as follows: if $b=(b_1,b_2,\hdots )$ and $b' = (b'_1,b'_2,\hdots )$ then
\[
(b.b')_n = \sum_{lcm(i,j)=n} (i,j) b_i b'_j \]
The interpretation map $itp$, which is a ringmorphism, sends $b=(b_1,b_2,\hdots )$ to the element
 of $\hat{B}(C)$ given by $\sum_{n=1}^{\infty} b_n [C_n]$ (where $C_n$ is the $C$-orbit of length $n$) and can thus be written as the difference  $[X_+]-[X_-]$ of two almost-finite $C$-sets, $X_+$ corresponding to the positive $b_n$ and $X_-$ to minus the negative $b_n$. The composition of the interpretation map with $\hat{\phi}$ is the {\em ghost map} \index{ghost map} $gh~:~b \mapsto d$ where $d_n = \phi_{C^n}(X(b)) = \phi_{C^n}(X_+(b)) - \phi_{C^n}(X_-(b))$ and the sequence of integers $d$ is related to that of $b$ via the formula
 \[
 \prod_{n=1}^{\infty} (\frac{1}{1-t^n})^{b_n} = exp(\int \sum_{n=1}^{\infty} d_n t^{n-1} dt) \]
 that is, $gh(b) = \hat{b}$ where $\hat{b}_n = \sum_{i| n} i b_i$.
  If $X$ is an almost-finite $C$-set then so is its $n$-th symmetric power
\[
S^n X = \{ g: X \rTo \N~|~g~\text{has finite support and}~\sum_{x \in X} g(x)=n \} \]
on which $C$ acts via $(c.g)(x) = g(c^{-1}.x)$. Dress and Siebeneicher prove in \cite{DressSieb89} that the map induced by
\[
s_t~:~\hat{B}(C) \rTo w(\Z) \qquad [X] \mapsto s_t(X)=1+\wis{fix}(S^1(X))t+\wis{fix}(S^2(X))t^2+ \hdots \]
is an isomorphism of rings and if $s_t(X(b)) = 1 + \sum_{n=1}^{\infty} a_n t^n = s_{\pmb{a}}(t)$ then we have that
\[
\prod_{n=1}^{\infty} (\frac{1}{1-t^n})^{b_n} = 1 + \sum_{n=1}^{\infty} a_n t^n \]
This allows us in the case of $w(\Z)$ to compute the product combinatorial. If $s_t(X(b)) = s_{\pmb{a}}(t)$ and $s_t(X(d))=s_{\pmb{c}}(t)$ then we have in $w(Z)$ that
\[
s_{\pmb{a}}(t) \otimes s_{\pmb{c}}(t) = s_t(X(b) \times X(d)) \]
If $m \in \N$ define the {\em congruence maps} \index{congruence maps}
\[
m^{(C)} = \{ g~:~C \rTo \{ 1,2,\hdots,m \}~|~\exists n \in \N~:~z_1.z_2^{-1} \in C^n = \langle c^n \rangle \Rightarrow g(z_1)=g(z_2) \} \]
Observe that $m^{(C)}$ is again an almost finite $C$-set under the action $c.g(z) = g(c^{-1}.z)$ and one verifies that $\phi_{C^n}(m^{(C)}) = m^n$. The map $m \mapsto m^{(C)}$ from $\N_0$ to $\hat{B}(C)$ extends to a map $\Z \rTo \hat{B}(C)$ and $\hat{\phi}(m^{(C)}) = (m,m^2,m^3,\hdots )$. We will now extend this map to a map $W(\Z) \rTo \hat{B}(C)$.
If $X$ is an almost finite $C$-set and $n \in \N$ then we define its {\em induction} \index{induction} with respect to the $n$-th power map $\sigma_n~:~C \rTo C$ given by $c \mapsto c^n$, $ind_n(X)$ as the set of $C$-orbits in $C \times X$ under the action $c.(c',x) = (c'c^{-n},c.x)$. Again, $ind_n(X)$ becomes an almost finite $C$-set via the action $c.\Oscr(c',x) = \Oscr(c.c',x)$ and one verifies that
\[
ind_n(X_1 \sqcup X_2) = ind_n(X_1) \sqcup ind_n(X_2) \quad \text{and} \quad ind_n(C_i) = C_{ni} \]
and have that $\phi_{C^m}(ind_n(X)) = n.\phi_{C^{m/n}}(X)$ if $n | m$ and is zero otherwise. This then gives a ring-isomorphism
\[
\tau~:~W(\Z) \rTo \hat{B}(C) \qquad q=(q_1,q_2,\hdots) \mapsto \sum_{n=1}^{\infty} ind_n(q_n^{(C)}) \]
If $\tau(q) = X(b)$ then the sequences of integers $q$ and $b$ are related via the formula
\[
\prod_{n=1}^{\infty} \frac{1}{1-q_n t^n} = \prod_{n=1}^{\infty} ( \frac{1}{1-t^n})^{b_n} \]

If $X$ is an almost-finite $C$-set, then $res_n(X)$ is the restriction to the subgroup $\langle c^n \rangle$. that is, $X=res_n(X)$ but with a new action $\odot$ defined by $c \odot x = c^n.x$. Clearly, $res_n$ is compatible with disjoint union and direct product and hence is an endomorphism
\[
res_n~:~\hat{B}(C) \rTo \hat{B}(C) \]
which are the {\em Adams operations} \index{Adams operations} on $\hat{B}(C)$ and this family of commuting endomorphisms of $\hat{B}(C)$ corresponds to the family of commuting endomorphisms $\Psi^n$ on $w(\Z)$. Similarly, the Verschiebungs additive maps on $w(\Z)$ are given by induction from the subgroup $\langle c^n \rangle$. 
Induction and restriction satisfy the following properties
\begin{itemize}
\item{$res_n(C_m) = (n,m) C_{[n,m]/n}$, where $[n,m]=lcm(n,m)$}
\item{$ind_n(C_m) = C_{nm}$}
\item{$\wis{ker}(res_n) = \{ x \in \hat{B}(C)~|~\phi_{C^m}(x) = 0~\forall n | m \}$}
\item{$\wis{im}(ind_n) = \{ x \in \hat{B}(C)~|~\phi_{C^m}(x)=0~\forall n \nmid m \}$}
\end{itemize}
Similarly, one can make Frobenius and Verschiebungs operators explicit on the necklace algebra $Nr(\Z)$. Define the Frobenius ringmorphisms $f_n~:~Nr(\Z) \rTo Nr(\Z)$ by
\[
f_n(b_1,b_2,\hdots ) = (\sum_{[n,i]=n} (n,i)b_i, \sum_{[n,i]=2n}(n,i) b_i, \hdots ) \]
and the Verschiebungs additive morphisms $v_n~:~Nr(\Z) \rTo Nr(\Z)$ via
\[
v_n(b_1,b_2, \hdots) = (\underbrace{0,\hdots,0}_{n-1},b_1,\underbrace{0,\hdots,0}_{n-1},b_2,\hdots ) \]
and these Frobenius and Verschiebungs operations $f_n$ and $v_n$ commute with induction and restriction maps $ind_n$ and $res_n$ on $\hat{B}(C)$.

In retrospect, the appearance of Burnside rings in $\F_1$-geometry is not surprising. Remember from the Smirnov-Kapranov paper \cite{SmirnovKapranov96} that $GL_n(\F_1) \simeq S_n$, so for any group $G$ an $n$-dimensional representation of $G$ over $\F_1$ would be a group-morphism $G \rTo S_n$, that is, a permutation representation of $G$, or equivalently, a finite $G$-set. If $G$ is an infinite discrete group, this says that any finite dimensional $\F_1$-representation of $G$ factors as a permutation representation though a finite group-quotient and hence determines an element in the Burnside-ring $B(\hat{G})$ of the profinite completion of $G$. In the special case when $G=C$ we can write $C$ additively (that is, $C=\Z$) and its $\C$-representations are of course all $1$-dimensional and parametrized by $\C^*$. The $\F_1$-representations are then the representations of the profinite completion $\hat{\Z}_+$ and its $\C$-points are precisely the roots of unity! Further, for completed Burnside rings we have $\hat{B}(G) = \hat{B}(\hat{G})$ so in our case $w(\Z) = \hat{B}(C) = \hat{B}(Z) = \hat{B}(\hat{\Z})$.

In \cite{DressSieb88} Andreas Dress and Christian Siebeneicher have extended the Witt construction to the profinite completion $\hat{G}$ of an arbitrary discrete group $G$ (and in fact even to arbitrary profinite groups). Let $cosg(\hat{G})$ be the set of conjugacy classes of open subgroups of $\hat{G}$ (that is, the conjugacy classes of subgroups of $G$ of finite index), then one can consider the covariant functor
\[
W_G~:~\wis{rings} \rTo \wis{rings} \qquad A \mapsto A^{cosg(\hat{G})} \]
and they show that with respect to this functor we have an isomorphism between $W_G(\Z)$ and the Burnside ring $\hat{B}(\hat{G})$ of almost finite $G$-sets. Moreover, the rings $W_G(R)$ all have Frobenius-like and Verschiebungs-like morphisms to (and from) $W_U(R)$, for any subgroup $U$ of $G$ of finite index. The Frobenius and Verschiebung maps
\[
W_G(R) \rTo^{\Psi^U} W_U(R) \quad \text{and} \quad W_U(R) \rTo^{V_U} W_G(R) \]
 are defined by restriction resp. induction. Clearly, in the case when $G=\Z$ all cofinite subgroups are isomorphic to $\Z$ giving rise to the Frobenius-lifts {\em endo}morphisms and corresponding Verschiebungs-operations on $w(R)$. 
 
 This raises the exciting prospect of extending or modifying Borger's $\lambda$-rings approach to $\F_1$-geometry to other categories $\wis{rings}_G$ of commutative rings with suitable morphisms to/from a collection of rings (for any conjugacy class of a cofinite subgroup of $G$) such that the Dress-Siebeneicher-Witt functor $W_G$ is a right adjoint functor to the forgetful functor $\wis{rings}_G \rTo \wis{rings}$. We expect such an approach to be fruitful when one starts with the braid group $B_3$ or its quotient the modular group $PSL_2(\Z)$, which may also clarify the role of Conway's big picture, which after all was intended to better understand cofinite subgroups of the modular group.
 
\section{What is $\wis{Spec}(\Z)$ over $\F_1$?}

So we can compute explicitly with $w(\Z)$ and know that $\wis{Spec}(\Z) / \F_1$ is the geometric object associated to $w(\Z)$, but what is this object and can we make sense of Smirnov's covering maps $\wis{Spec}(\Z) \rOnto \PP^1_{\F_1}$? 

We have a candidate for the geometric object, namely the $\lambda$-spectrum of $w(\Z)$
\[
\wis{Spec}(\Z)_{\F_1} = \wis{Spec}_{\lambda}(w(\Z)) = \{ \wis{ker}(w(\Z) \rTo A)~|~A~\text{reduced $\lambda$-ring}~\} \]
If $\wis{Spec}(\Z)$ would behave as a 'curve' over $\F_1$, one would expect the $\lambda$-spectrum to contain many geometric points over $\overline{\F_1}$. However, we will soon see that
\[
\wis{max}_{\lambda}(w(\Z)) = \emptyset = \wis{max}_{cycl}(w(\Z)) \]
In fact a similar result holds for any $w(R)$. 

The fact we will use is that the Verschiebungs-operators survive $\lambda$-ring quotients $A=w(R)/I$ which have no additive torsion. I thank Jim Borger for communicating this to me.
Clearly, there are additive maps 
\[
v_n~:~w(R) \rTo^{V_n} w(R) \rOnto A \]
and we have to show that $\wis{ker}(v_n) \subset I$. Because the $\Psi^p$ are lifts of the Frobenius, there is a unique map $d$ on $w(R)$ such that for all $s(t) \in w(R)$ we have the identity
\[
(s(t))^{\otimes p} + [p] d(s(t)) = \Psi^p(s(t)) \]
and hence any $\lambda$-ideal $I$ must be preserved by $d$. Assume $s(t) \in \wis{ker}(v_{np})$, that is $V_{np}(s(t)) \in I$, then from the identities
\[
(V_{np}(s(t))^{\otimes p} + [p] d(V_{np}(s(t))) = \Psi^p(V_{np}(s(t))) = \Psi^p \circ V_p (V_n(s(t))) = [p] V_n(s(t)) \]
it follows that the left-hand side is contained in $I$ so must be the right-hand side and because $A=w(R)/I$ has no additive torsion it follows that $V_n(s(t)) \in I$, so $v_n(s(t)) \in \wis{ker}(v_n)$. As we can repeat this process for any prime factor $p$ of $m=np$ it follows that if $s(p) \in \wis{ker}(v_m)$ then $s(t) \in \wis{ker}(v_1)=I$. 
Thus, if $A$ is a $\lambda$-ring quotient of $w(R)$ without additive torsion, $A$ is equipped not only with ring-endomorphisms $\Psi^n$ but also with additive morphisms $v_n$ satisfying all the properties the Frobenius and Verschiebungs operators satisfy on $w(R)$, indicating that $A$ must itself be close to a witty ring.

Now assume that $A$ is \'etale over $\F_1$ and hence of finite rank over $\Z$. Recall from \cite{LenstraWitt} that we can also define the ringstructure of $w(R)$ as the inverse limit 
\[
w(R) = \underset{\leftarrow}{lim}~w_n(R) \quad \text{with} \quad  w_n(R)=\wis{ker}((R[t]/(t^{n+1})^* \rTo R^*) \]
As $A$ is finite over $\Z$. the ringmorphism $w(R) \rOnto A$ factors through a $w_n(R)$ for some $n \in \N$. But this means that $V_n(w(R))$ is contained in the ideal $I$, in particular $v_n(1) \in I$ but then from the argument given before we conclude that $1 \in I$ and hence that $A=0$. That is, witty rings $w(R)$ do not have torsion-free $\lambda$-ring quotients, finite over $\Z$.

That is, $\lambda$-spectra of witty rings do not have geometric points and hence behave very unlike $\F_1$-geometric objects of finite dimension. Still, the $\lambda$-spectrum has many other points, in fact we can identify the usual prime spectrum $\wis{Spec}(R)$ with a subset of witty points in $\wis{Spec}_{\lambda}(w(R))$
\[
\wis{Spec}(R) \simeq \wis{Spec}_w(w(R)) = \{ \wis{ker}(w(R) \rTo w(Q(R/p))~|~p \in \wis{Spec}(R) \} \]
where $Q(R/p)$ denotes the field of fractions of the domain $R/p$. 

Let us work out what the witty ring $w(F)$ of a field $F$ is. If $\overline{F}$ is algebraically closed, then by construction we have an inclusion of multiplicative groups
$\overline{F}^* \rInto w(\overline{F})^*$ determined by $a \mapsto \frac{1}{1-at}$
which extends to a ringmorphism on the group-algebra of $\overline{F}$, 
$\Z[\overline{F}^*] \rInto^L w(\overline{F})$ 
with image the set of all rational formal power series $\frac{\prod_i (1-\alpha_it)^{e_i}}{\prod_j (1- \beta_j t)^{f_j}} = L(\sum_j f_j \beta_j - \sum_i e_i \alpha_i)$. In other words, we have a suitably dense subring of $w(\overline{F})$ isomorphic to the integral group-algebra $\Z[\overline{F}^*]$. The absolute Galois group $G=Gal(\overline{F}/F)$ acts on both rings giving an inclusion of rings
\[
\Z[\overline{F}^*]^G = Div(\PP^1_F - \{ 0,\infty \}) \rInto w(F) \]
where $Div(\PP^1_F - \{ 0,\infty \})$ is the Abelian group of divisors on $\PP^1_F - \{ 0,\infty \}$, that is all formal finite combinations $\sum n_i [f_i]$ with $n_i \in \Z$ and the $f_i$ irreducible monic polynomials in $F[x,x^{-1}]$ and it gets an induced multiplication (and even $\lambda$-ring structure) from the group-algebra structure $\Z[\overline{F}^*]$. The Frobenius lifts on $\Z[\overline{F}^*]$ are the linearizations of the multiplicative group endomorphisms sending $a \mapsto a^n$ for all $a \in \overline{F}^*$. 

In the special case of $\F_p$ we have seen before that we can identify the multiplicative group $\overline{\F}_p^*$ with the group of all roots of unity $\pmb{\mu}^{(p)}$ of order prime to $p$ and hence we get a dense subring
\[
\Z[\pmb{\mu}^{(p)}]^{\hat{\Z}_+} = Div(\PP^1_{\F_p} - \{ 0,\infty \}) \rInto w(\F_p) \]
and hence a surprise guest re-appearance of the fiber $\PP^1_{\F_p}$ of the structural map $\PP^1_{\Z} \rOnto \wis{Spec}(\Z)$ in the description of the witty point in $\wis{Spec}_{\lambda}(w(\Z))$ determined by the $\lambda$-ring morphism $w(\Z) \rOnto w(\F_p)$, somewhat closing the circle of thoughts we began by looking at Mumford's drawings!

Still, there's the eternal problem of finding a natural identification between $\overline{\F}_p^*$ and $\pmb{\mu}^{(p)}$. We will briefly sketch how this can be done 'in principle' using ordinal numbers. In \cite{ConwayONAG} John H. Conway identified the algebraic closure of $\F_2$ with the set of all ordinal numbers smaller than $\omega^{\omega^{\omega}}$ equipped with nim-addition and multiplication. Later Joseph DiMuro extended this to identify the algebraic closure of $\F_p$ with $\omega^{\omega^{\omega}}$ in \cite{DiMuro11}. We will recall the characteristic two case and refer to \cite{DiMuro11} for the general case. 

To distinguish the nim-rules from addition and multiplication of ordinal numbers, we will denote the later ones enclosed in brackets. So, for example $[\omega^2]$ will be the ordinal number, whereas $\omega^2$ will be the square of the ordinal number $[\omega]$ in nim-arithmetic. These nim-rules can be defined on all ordinals as follows
\[
\alpha+\beta=\wis{mex}(\alpha'+\beta,\alpha+\beta') \quad \text{and} \quad \alpha.\beta=\wis{mex}(\alpha'.\beta + \alpha.\beta'-\alpha'.\beta') \]
where $\alpha'$ (resp. $\beta'$) ranges over all ordinals less than $\alpha$ (resp. than $\beta$) and $\wis{mex}$ stands for 'minimal excludent' of the given set, that is the smallest ordinal not contained in the set. Observe that these definitions allow us compute with ordinals inductively. Computing the sum of two ordinals is easy: write each one uniquely as sum of ordinal numbers $\alpha=[2^{\alpha_0} + 2^{\alpha_1} + \hdots + 2^{\alpha_k}]$ then to compute $\alpha+\beta$ we delete powers appearing in each factor and take the Cantor ordinal sum of the remaining sum (for finite ordinals this is the common nim-addition 'adding binary expressions without carry'). To compute multiplication of ordinals, introduce the following special elements
\[
\kappa_{2^n} = [2^{2^{n-1}}] \quad \text{and for primes $p>2$} \quad \kappa_{p^n} = [2^{\omega^k.p^{n-1}}] \]
where $k$ is the number of primes strictly smaller than $p$. Because $[2^{\alpha_0}+\hdots+2^{\alpha_k}]=[2^{\alpha_0}]+\hdots+[2^{\alpha_k}]$ we can multiply two ordinals $< [\omega^{\omega^{\omega}}]$ if we know how to compute products $[2^{\alpha}].[2^{\beta}]$ with $\alpha,\beta < [\omega^\omega]$. Each such $\alpha,\beta$ can be expressed uniquely as
\[
[\omega^t.n_t + \omega^{t-1}.n_{t-1} + \hdots + \omega.n_1 + n_0] \]
with $t$ and all $n_k$ natural numbers. If we write $n_k$ is base $p$ where $p$ is the $k+1$-st prime number, that is $n_k = [\sum_j p^j.m(j,k)]$ for $0 \leq m(j,k) < p$ then we can write any $2$-power smaller than $[\omega^{\omega^{\omega}}]$ as a decreasing finite product $[\prod_q \kappa_q^{m(q)}]$ with $0 \leq m(q) < p$ and $q$ a power of $p$. Conway has shown that for these $[\prod_q \kappa_q^{m(q)}] = \prod_q \kappa_q^{m(q)}$ which allows us to compute all products except when $[m(q)+m'(q)] \geq p$ for some $q$ appearing. Thus it remains to specify the ordinals $\kappa_q)^p$ and here Conway proved the rules
\[
(\kappa_{2^n})^2 = \kappa_{2^n} + \prod_{1 \leq i \leq n} \kappa_{2^i} \quad (\kappa_p)^p= \alpha_p \quad \text{and} \quad (\kappa_{p^n})^p = \kappa_{p^{n-1}} \]
for $p$ an odd prime and $n \geq 2$. Conway calculated the first few $\alpha_p$, for example $\alpha_3=2,\alpha_5=4, \alpha_7=[2^{\omega}]+1$ etc. and then Hendrik Lenstra \cite{Lenstra77} gave an explicit algorithm to compute the $\alpha_p$ and managed to determine then for all $p \leq 43$. Today we know all $\alpha_p$ for $p \leq 293$ with only a few exceptions. In principle this allows us to determine the ordinal number corresponding to any realistic occurring element in $\overline{\F}_2$. Similarly, DiMuro proved that $\overline{\F}_p$ can be identified with $[\omega^{\omega^{\omega}}]$ and listed the values for the $\alpha_q$ in those cases for primes $q \leq 43$ and $p \leq 11$. 

Using this correspondence we can now construct a one-to-one correspondence $\overline{\F}_p^* \leftrightarrow \pmb{\mu}^{(p)}$, which we will illustrate in the case $p=2$. Conway showed that the ordinals $[2^{2^n}]$ form a subfield isomorphic to $\F_{2^{2^n}}$ and so there is a consistent embedding of the quadratic closure of $\F_2$ into roots of unity by starting with $[2]$ being the smallest ordinal generating the multiplicative group of the subfield $[2^2]$ (of order $3$) and taking it to be $e^{2 \pi i}{3}$, for the next subfield $[2^{16}]$ we have to look for the smallest ordinal $[k]$ such that $[k]^5 = 2$, which turns out to be $[4]$ which then corresponds to $e^{2 \pi i/15}$ and the correspondence between $\F_{2^4} \leftrightarrow \pmb{\mu}_{15}$ is depicted below (together with the addition and multiplication tables of $[16]$ to verify the claims).
\[
\includegraphics[height=6.5cm]{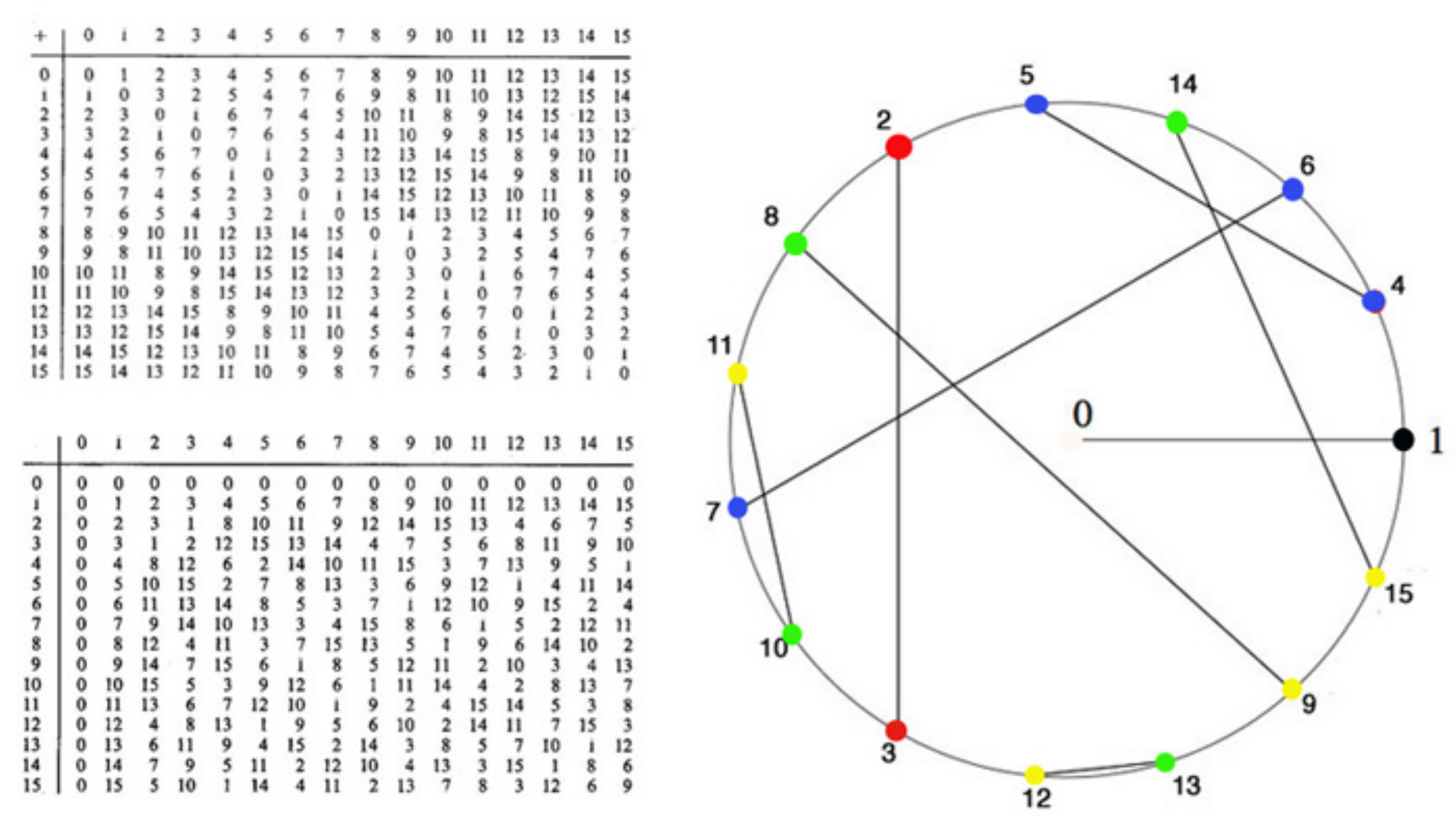}
\]
and we have indicated the different orbits under the Frobenius $x \mapsto x^2$ with different colors. There are two orbits of size one : $\{ 0 \}$ corresponding to $x$ and $\{ 1 \}$ corresponding to $x+1$. One orbit $\{ 2,3 \}$ of size two corresponding to the irreducible polynomial $(x-2)(x-3) = x^2+x+1$
and three orbits of size four corresponding to the three irreducible monic polynomials in $\mathbb{F}_2[x]$ of degree four, for example $\{ 4,6,5,7 \} \leftrightarrow x^4+x+1$. Iterating this procedure we get an explicit embedding of the quadratic closure of $\F_2$ into roots of unity (the relevant generators for the next stages are $32$ then $1051$ 
then $1361923$ and $1127700028470$). Having obtained an explicit identification of the quadratic closure of $\F_2$ in the roots of unity we can then proceed by associating to $[\omega]$ to root $e^{2 \pi i/9}$ as $[\omega]^3=[2]$, mapping $[\omega^{\omega}]$ to $e^{2 \pi i/75}$ as $[\omega^{\omega}]^5 = [4]$ and so on until we have identified $\overline{\F}_2^*$ with $\pmb{\mu}^{(2)}$. This then allows us to associate to a schematic point of $\mathbb{A}^1_{\F_2}$, that is an irreducible monic polynomial in $\F_2[x]$ the root of unity corresponding to the smallest ordinal in the Frobenius-orbit associated to the polynomial. So, for example, to $x^4+x^3+x^2+x+1$ one assigns $e^{12 \pi i/15}$ as the roots of the polynomial are the ordinals $\{ [8],[10],[13],[14] \}$. Once again, one can repeat these arguments for the algebraic closures $\overline{\F}_p$ using the results from \cite{DiMuro11}.

\section{What is the map from $\overline{\wis{Spec}(\Z)}$ to $\PP^1_{\F_1}$?}

In the foregoing sections we have recalled some of the successes of Borger's approach to absolute geometry via $\lambda$-rings. For example, the identification of the \'etale site of $\F_1$ with the category of finite sets equipped with an action of the monoid $\hat{\Z}_{\times}$ is one of the most convincing theories around vindicating Smirnov's proposal that one should interpret $\pmb{\mu}_{\infty}$ as the algebraic closure of $\F_1$. Further, with this $\lambda$-ring approach we obtain roughly the same class of examples provided by all other approaches to $\F_1$-geometry, such as affine- and projective spaces, Grassmannians, toric varieties, among others. In addition, we can associate a space of geometric points as well as a new topology to such an $\F_1$-geometric object $X$. For, assume that $X$ is locally controlled by a $\lambda$-ring $R$, then locally its {\em geometric points} \index{geometric point}  correspond to kernels of $\lambda$-ring morphisms $R \rOnto S$ where $S$ is \'etale over $\F_1$, among which are the {\em cyclotomic coordinates} \index{cyclotomic coordinates} which are the special points obtained by taking $S = \Z[x]/(x^n-1) = \Z[\pmb{\mu}_n]$. But, as $S$ is finite projective over $\Z$, these kernels are {\em not} maximal ideals of the $\lambda$-ring $R$, but rather sub-maximal ones, entailing that two such kernels no longer have to be co-maximal. This then leads to a adjacency- (or clique-) relation among the corresponding geometric points which gives us the {\em Habiro topology} \index{Habiro topology} on $\wis{max}_{\lambda}(R)$. This new topological feature encodes the fact that the closed subschemes of the usual integral affine scheme $\wis{Spec}(R)$ corresponding to the kernels of the two geometric points intersect over certain prime numbers $p$. As an example, we have seen that the cyclotomic points of $\PP^1_{\Z}$ (for the toric $\lambda$-structure) give us indeed the proposal that
\[
\PP^1_{\F_1} = \{ [0],\infty] \} \cup \{ [n]~:~n \in \N_0 \} \]
where two cyclotomic points $[n]$ and $[m]$ are adjacent if and only if their quotient is a pure prime power, leading to the Habiro topology on the roots of unity $\pmb{\mu}_{\infty}$. Further, the $\lambda$-ring structure, that is the commuting family of endomorphisms $\{ \Psi^n~:~n \in \N_0 \}$ can be viewed as descent data from $\wis{Spec}(\Z)$ to $\wis{Spec}(\F_1)$, and hence conversely we can view the process of forgetting the $\lambda$-ring structure as the base-extension functor $- \times_{\wis{Spec}(\F_1)} \wis{Spec}(\Z)$. In particular we can now make sense of the identity
\[
\PP^1_{\F_1} \times_{\wis{spec}(\F_1)} \wis{Spec}(\Z) = \PP^1_{\Z} \]
where the right-hand side is the usual integral scheme $\PP^1_{\Z}$, without emphasis on the toric $\lambda$-structure.

But Borger's proposal really shines in that it allows us to make sense of what the Weil restriction to $\wis{Spec}(\F_1)$ is of {\em any} integral scheme. Indeed, the witty functor $w~:~\wis{rings} \rTo \wis{rings}_{\lambda}$ is the right adjoint of the forgetful functor (which we have seen is base-extension) and hence if the integral scheme $X$ is locally of the form $\wis{Spec}(R)$, then $X/\F_1$ is locally the geometric object corresponding to the $\lambda$-ring $w(R)$. However, such rings do not have geometric points as before, so we have a dichotomy among the geometric $\F_1$-objects which resembles that in noncommutative algebraic geometry between algebras having plenty of finite dimensional representations, and algebras that have no such representations. Geometric $\F_1$-objects are either the restricted class of combinatorial controlled integral schemes allowing a $\lambda$-structure, or the class of infinite dimensional objects corresponding to Witt schemes of integral schemes. Still, the ordinary integral scheme structure survives this witty-fication as $\wis{Spec}(R)$ can be embedded in the $\lambda$-spectrum $\wis{Spec}_{\lambda}(w(R))$ via the kernels of the $\lambda$-ring maps $w(R) \rTo w(R/P) \rTo w(Q(R/P))$ for any prime ideal $P$ of $R$. As an example, $\wis{Spec}(\Z)$ is the $\F_1$-geometric object corresponding to the Burnside ring $w(\Z) = \hat{B}(C)$ which does indeed contain the proposal that
\[
\wis{Spec}(\Z)/\F_1 = \{ (p)~:~p~\text{a prime number} \} \]
where we view the prime number $p$ as corresponding to the $\lambda$-ring morphism
\[
w(\Z) \rTo w(\F_p) \approx \wis{Div}(\PP^1_{\F_p} - \{ [0],[\infty] \}) \]
Although these two classes of geometric $\F_1$-objects are quite different, we can still make sense of morphisms between them, as they have to be locally given by $\lambda$-ring morphisms. In particular,
let us investigate whether we can make sense of Smirnov's maps 
\[ q= \frac{a}{b}~:~\wis{Spec}(\Z) \rTo \PP^1_{\F_1} \]
 in Borger's $\lambda$-rings approach to $\F_1$-geometry, that is, this map should locally be determined by a $\lambda$-ring morphism. With $\PP^1_{\F_1}$ we mean the cyclotomic points of the integral scheme $\PP^1_{\Z}$ equipped with the toric $\lambda$-ring structure. Because $(a,b)=1$, we can cover $\PP^1_{\Z}$ with the prime spectra of two $\lambda$-rings
 $\PP^1_{\Z} = \wis{Spec}(\Z_b[x]) \cup \wis{Spec}(\Z_a[x^{-1}])$, and therefore
 \[
 \PP^1_{\F_1} = \wis{Spec}_{cycl}(\Z_b[x]) \cup \wis{Spec}_{cycl}(\Z_a[x^{-1}]) \]
 Further, we have seen that $\wis{Spec}(\Z) / \F_1$ should be viewed as $\wis{Spec}_{\lambda}(w(\Z))$ which we can cover as $\wis{Spec}_{\lambda}(w(\Z_b)) \cup \wis{Spec}_{\lambda}(w(\Z_a))$. Now, consider the $\lambda$-ring morphisms
 \[
 \begin{cases} \Z_b[x] \rTo w(\Z_b) &~x \mapsto \frac{1}{1-\frac{a}{b}t} \\
 \Z_a[x^{-1}] \rTo w(\Z_a) \qquad &~x^{-1} \mapsto \frac{1}{1-\frac{b}{a}t }
 \end{cases}
 \]
 which coincide on the intersection with the $\lambda$-morphism
$ \Z_{ab}[x,x^{-1}] \rTo w(\Z_{ab})$ determined by $x \mapsto \frac{1}{1-\frac{a}{b}t}$. So, in order to investigate the associated geometric map
\[
\wis{Spec}(\Z_b) \simeq \wis{Spec}_w(w(\Z_b)) \rTo \wis{Spec}_{cycl}(\Z_b[x]) \]
we have to look for any prime $p$ not dividing $b$ to the composition $\Z_b[x] \rTo w(\F_p)$ which maps $x$ to $\tfrac{1}{1-\frac{\overline{a}}{\overline{b}}t}$ and hence for every $n \in \N_0$
\[
x^n \mapsto \frac{1}{1-(\frac{\overline{a}}{\overline{b}})^nt} \]
If $\overline{a}/\overline{b}$ has order $n$ in $\F_p$ this says that $x^n$ maps to $1/(1-t) = \boxed{1} \in w(\F_p)$ and if $p | a$ then $x$ is mapped to $\tfrac{1}{1-0t}=1=\boxed{0} \in w(\F_p)$. Further, if $p | b$ we get in the composition $\Z_a[x^{-1}] \rTo w(\F_p)$ that $x^{-1}$ is mapped to $\frac{1}{1-0t} = 1 = \boxed{0} \in w(\F_p)$. So, if we write $[p]$ for the witty-point corresponding to the kernel of $w(\Z) \rOnto w(\F_p)$ we get indeed that
\[
[p] \mapsto \begin{cases} [0]~\text{if $p|a$} \\
[\infty]~\text{if $p|b$} \\
[n]~\text{if $n$ is minimal such that $a^n-b^n \equiv 1~\wis{mod}~p$}
\end{cases} \]
which indeed coincides with Smirnov's proposal.

\vskip 4mm

\noindent
{\bf Acknowledgements : } I thank Jim Borger and Jack Morava for several illuminating emails. These notes are based on a rather chaotic master course given in Antwerp in 2011/12. I thank the students for their patience and inspiring enthusiasm and, in particular, Pieter Belmans for pointing me to Mumford's $\PP^1_{\Z}$-picture as well as for generous help with Sage/TikZ in order to produce some of the pictures.

\printindex
\end{document}